%% file: main.tex
\documentclass[a4paper,11pt]{article}

\input{packages}%
\input{commands}%

\begin{document}

\title{A Framework for Directional and Higher-Order \\ Reconstruction in Photoacoustic Tomography}
\author[1,2]{Yoeri~E.~Boink\footnote{\href{mailto:y.e.boink@utwente.nl}{y.e.boink@utwente.nl}}}
\author[3]{Marinus~J.~Lagerwerf}
\author[2]{Wiendelt~Steenbergen}
\author[1]{Stephan~A.~van~Gils}
\author[2]{Srirang~Manohar}
\author[1]{Christoph~Brune}
\affil[1]{Department of Applied Mathematics, University of Twente, The Netherlands}%
\affil[2]{Biomedical Photonic Imaging Group, University of Twente, The Netherlands}%
\affil[3]{Computational Imaging Group, Centrum Wiskunde \& Informatica, The Netherlands}

\date{\today}

\maketitle

\begin{abstract}
\noindent Photoacoustic tomography is a hybrid imaging technique that combines high optical tissue contrast with high ultrasound resolution. Direct reconstruction methods such as filtered backprojection, time reversal and least squares suffer from curved line artefacts and blurring, especially in case of limited angles or strong noise. In recent years, there has been great interest in regularised iterative methods. These methods employ prior knowledge on the image to provide higher quality reconstructions. However, easy comparisons between regularisers and their properties are limited, since many tomography implementations heavily rely on the specific regulariser chosen. To overcome this bottleneck, we present a modular reconstruction framework for photoacoustic tomography. It enables easy comparisons between regularisers with different properties, e.g. nonlinear, higher-order or directional. We solve the underlying minimisation problem with an efficient first-order primal-dual algorithm. Convergence rates are optimised by choosing an operator dependent preconditioning strategy. Our reconstruction methods are tested on challenging 2D synthetic and experimental data sets. They outperform direct reconstruction approaches for strong noise levels and limited angle measurements, offering immediate benefits in terms of acquisition time and quality. This work provides a basic platform for the investigation of future advanced regularisation methods in photoacoustic tomography.\\
\\
\textbf{keywords:} photoacoustic tomography, variational image reconstruction, total generalised variation, directional regularisation, compressive sampling, convex optimisation.
\end{abstract}

\newpage
\section{Introduction}
Photoacoustic tomography (PAT), also known as optoacoustic tomography, is an intrinsically hybrid imaging technique that combines the high spectroscopic contrast of tissue constituents to light, with the high resolution of ultrasound imaging techniques. Tissue is illuminated with nanosecond laser pulses, causing some of the optical energy to be absorbed and converted into heat. This leads to thermoelastic expansion and creates ultrasound waves that are detected at the boundary. These detected signals are employed to reconstruct the spatial distribution of absorbed optical energy inside the tissue. Common methods applied are filtered backprojection (FBP) \cite{Kruger1995, Haltmeier2005, Xu2005, Willemink2010} or time reversal \cite{Hristova2008,Treeby2010}. In this work, we develop a general variational framework that incorporates prior knowledge on the reconstruction, offering higher quality reconstructions than direct methods and providing robustness against noise and compressive sampling.

The optical absorption coefficient in tissue varies due to the relative presence of various components such as haemoglobin, oxy-haemoglobin, melanin, lipids and water. These chromophores have specific spectral signatures, and the use of multi-wavelength PAT can potentially reveal molecular specific information. This is directly linked to tissue physiopathology and function and has diagnostic potential. Photoacoustics is being researched for applications in various fields of biomedicine \cite{Zhou2016} such as brain imaging in small animals \cite{Dean-Ben2017}, breast cancer imaging in humans \cite{Heijblom2015, Toi2017}, and imaging inflamed human synovial joints in rheumatoid arthritis \cite{Chamberland2010,Es2014}. In most of these applications, detection of signals from haemoglobin and oxy-haemoglobin enables the visualisation of blood vessels. Several disorders are characterised by a dysregulation of blood vessel function, but also with uncontrolled creation of blood vessels or angiogenesis.

The standard methods of filtered backprojection (FBP) and time reversal suffer from curved line artefacts, especially when the noise level is high or the placement of detectors for measurements (sampling) is coarse \cite{Rosenthal2013}. Recently, there has been intense interest in solving the PAT reconstruction problem iteratively with a specific focus on regularised reconstruction. The regularisers used for these reconstructions have many different forms, depending on the prior assumptions on the image. Total variation (TV) regularisation is a powerful tool for noise removal \cite{Rudin1992}, and generates the desired images with sharp edges. In \cite{Arridge2016b}, TV-regularised reconstructions on simulated data are presented, where an analytical adjoint operator for the k-space forward operator \cite{Treeby2010} is derived to model wave propagation. In \cite{Arridge2016a} a similar reconstruction method is applied to compressively sampled simulated and experimental data. Although from the theory it is not obvious \cite{Provost2009}, the TV basis appears to work well in the area of compressive sampling for the reconstruction of an absorbed energy density map with abrupt edges \cite{Guo2010}. Other works that combine total variation with compressive sampling in PAT can be found in \cite{Zhang2012} and \cite{Meng2012}.

In case of a heterogeneous fluence rate, such as an exponentially decaying one, a piecewise constant prior on the reconstruction is not realistic. Higher-order total variation methods, such as total generalised variation (TGV) \cite{Bredies2010} are better suited to deal with such data. The use of TGV shows great promise for image reconstruction in tomographic settings, such as MRI \cite{Knoll2011}, CT \cite{Yang2010, Niu2014}, ultrasound waveform tomography \cite{Lin2016} and diffusion tensor imaging \cite{Valkonen2013a}. Image denoising with TGV as a post-processing step has made its way into optical coherence tomography \cite{Duan2016} and optical diffraction tomography \cite{Krauze2016}. However, solving these optical reconstruction problems with TGV regularisation within the reconstruction algorithm is still an open research topic. A recent report \cite{Hammernik2017} shows TGV-regularised reconstructions in PAT on vascular data. Here, measurements are taken by making snapshots of the acoustic waves using a CCD camera.

Vascular structures show a strong anisotropy, that can be promoted by using directional wavelets \cite{Kingsbury1998} or curvelets \cite{Candes1999} as `building blocks' for the reconstruction. Provost and Lesage \cite{Provost2009} have shown that wavelets and curvelets are very sparse representations of the measurement data, and can help to recover anisotropic features in the reconstruction. This has been confirmed in \cite{Guo2010} through new phantom and in vivo experiments. 

To our knowledge, a general reconstruction framework for PAT, in which regularisers can be easily exchanged, does not exist. In this work, we give such a variational framework. Its applicability is demonstrated by using it for three reconstruction models, which contain TV, TGV and directional wavelet regularisers. While others based their methods on specific measurement geometries, our framework makes no assumption on the geometry, requires no specific physics modelling and can be applied to the PAT reconstruction problem in a 2- or 3-dimensional setting. We develop not only a TV, but also a TGV regularised reconstruction model, which can deal with the more realistic assumption of a decaying heterogeneous fluence rate. A reconstruction model using directional wavelet regularisation is developed for images containing anisotropic vascular structures. The numerical implementation with a primal-dual algorithm is a modular one: the regularisation or the data-fidelity can be chosen differently without having to change the structure. The algorithm has a convergence rate that is at least linear and that can be accelerated \cite{Pock2011} for specific choices of the regulariser. The performance of the framework is demonstrated on synthetic data, as well as on experimental data. Robustness against noise and compressive sampling is shown. 

The remainder of this paper is organised as follows: in section \ref{sec:var_meth} the variational method for the regularised reconstruction of PAT is derived. After writing the reconstruction problem as a saddle-point problem, the numerical implementation with the first-order primal-dual algorithm PDHGM is explained in section \ref{ch:4}. In section \ref{sec:exp_setup}, the experimental setup that is employed to generate our results is explained. Moreover, the specific forward model that is used for obtaining the synthetic data is derived. The digital and experimental phantoms that are created as test cases for our reconstruction framework are shown in section \ref{sec:phantoms}. After an extensive analysis of the method on challenging synthetic data in section \ref{sec:results} and experimental data in section \ref{sec:results_exp}, we conclude and discuss the results in section \ref{sec:discussion}.

\section{Variational methods using regularisation} \label{sec:var_meth}
We consider the image reconstruction problem of finding an estimate for the absorbed energy density map $u\in L^2(\Omega)$ from given (preprocessed) measurements $f\in L^2(\Sigma)$. Here $\Omega\subset\R^d$, where $d$ is the dimension of the space in which $u$ must lie and $\Sigma = S\times(0,T)$ is the space-time domain of the detection surface $S$ between times $0$ and $T$. Mathematically, the inverse problem can be formulated as the solution to the equation 
\begin{equation}\label{eq:objective}
Ku = f,
\end{equation}
where $K: L^2(\Omega) \mapsto L^2(\Sigma)$ is a bounded linear operator. This operator can be described in many ways, since there are multiple methods that model the forward process in photoacoustic imaging. For instance, it can be described as a k-space method that models linear acoustic propagation in time \cite{Treeby2010}. In our simulations and experiments, we use the spherical mean operator \cite{Kruger1995} to model the forward process. 

Instead of using a direct reconstruction method such as FBP, we solve the inverse problem in a variational setting. We consider the following minimisation problem:
\begin{align}\label{eq:min_prob}
\min_{u\in L^2(\Omega)}&\left\{\frac12\norm{Ku-f}_{L^2(\Sigma)}^2 +\alpha R(u) \right\}.
\end{align}
The first term is the $L^2$ data fidelity term, which ensures that the reconstructed image, after applying the forward operator $K$, is close to the noisy data $f$. The $L^2$-norm has been chosen because the photoacoustic measurements are predominantly affected by noise from the system electronics and thermal noise from the transducers, which are known to be of additive Gaussian type \cite{Zhang2012}. The second term is an $\alpha$ weighted regularisation term which becomes large when our prior assumptions on $u$ are violated. Note that this is a very general framework: as long as $R(u)$ is convex, we can put it in our numerical framework, which is explained in section \ref{ch:4}. 

Note that if we do not have any prior assumption on the data and we choose $R(u)=0$, we obtain for \eqref{eq:min_prob} the least squares (LS) reconstruction. The LS reconstruction can be evaluated by using a gradient descent algorithm or, if $K$ is a small enough matrix operator, by solving $u = (K^T K)^{-1}K^T f$ directly. However, since LS is very sensitive to noise in the data $f$, Tikhonov regularised least squares (LS-T) is often used: the term $R(u)$ in \eqref{eq:min_prob} then takes the form $\frac12 \norm{u}^2_{L^2(\Omega)}$. This can again be evaluated directly by solving $u = (K^T K + \alpha)^{-1}K^T f$. In section \ref{sec:results}, our models will be compared with LS-T. 

\subsection{Regularisation with total (generalised) variation}
In this work we use the concept of total generalised variation (TGV), as introduced by Bredies \textit{et al.} \cite{Bredies2010}. TGV is a generalisation of TV, in which higher order derivatives are considered instead of only first order derivatives in TV. Bredies \textit{et al.} proposed the following definition:
\begin{align}\label{eq:def_TV}
\TGV_\beta^k(u) = \sup_q\bigg\{\int_\Omega u~\div^k(q)\diff x~\Big|~q\in \C_0^k(\Omega,\text{Sym}^k(\R^d)),\norm{\div^l(q)}_\infty\leq\beta_l,~l=0,\hdots,k-1\bigg\}.
\end{align}
With the choice of $k$, the desired order of regularity can be obtained, while the parameter $\beta$ gives a weight on every regularity level. By choosing $k=1$ and $\beta=1$, the definition of TGV coincides with the definition of TV
\begin{align}\label{eq:def_TGV}
\TV(u) =&\sup_q\bigg\{\int_\Omega u~\div(q)\diff x~\Big|~q\in \C_0^1(\Omega,\R^d),~\norm{q}_\infty\leq 1\bigg\}.
\end{align}
In this work, we consider $\TV(u)$ and $\TGV^2_\beta(u)$ as regularisers. In our applications, we work with a discrete image $u$, on which a gradient or derivative is always numerically defined. Therefore, if we assume $u\in\C^1(\Omega)$ or $u\in\C^2(\Omega)$, the expressions for $\TV(u)$ and $\TGV_\beta^2(u)$ respectively can be simplified. After substituting these expressions in \eqref{eq:min_prob}, we obtain the following minimisation problems:
\begin{align}
&u_{\scriptscriptstyle \TV} = \argmin_{u\in\BV}\left\{\frac12\norm{Ku-f}_{L^2(\Sigma)}^2 +\alpha \norm{\nabla u}_{(L^1(\Omega))^d} \right\} \label{eq:min_TV},\\
&u_{\scriptscriptstyle \TGV} = \argmin_{\substack{u,v\in\BV}}\bigg\{\frac12\norm{Ku-f}_{L^2(\Sigma)}^2+\alpha \big\{ \norm{\nabla u - v}_{(L^1(\Omega))^d} +\beta\norm{\mathcal{E}(v)}_{(L^1(\Omega))^{d\times d}} \big\}\bigg\}, \label{eq:min_TGV}
\end{align} where $\BV = \{u\in L^2(\Omega)~|~\TV(u)<\infty\}$ and $\mathcal{E}$ is the symmetrised gradient \cite{Bredies2010}. In the TV-regularised functional \eqref{eq:min_TV}, the parameter $\alpha$ determines the smoothness of the solution. In the TGV-regularised functional \eqref{eq:min_TGV} one can choose the influence of first order and second order smoothness by choosing combinations of $\alpha$ and $\beta$. The effect of different combinations of $\alpha$ and $\beta$ on the reconstruction of one-dimensional TV and TGV eigenfunctions was analysed in \cite{Benning2013}. From now on, when TGV is mentioned, we mean the second order method $\TGV_\beta^2$.

\subsection{Regularisation with directional wavelets}
Blood vessel visualisation is an important aspect for the application of photoacoustic imaging in biomedicine. An image containing vascular structures can be sparsely represented by only a small number of elements in some basis or in a so called dictionary. There is a big amount of possibilities to represent this anisotropic data. A rough distinction can be made \cite{Rubinstein2010} between analytic transforms, such as the Gabor and wavelet transform; analytic dictionaries, such as curvelets and contourlets; and trained dictionaries, such as sparse PCA and K-SVD. 

In this work, we make use of an enhancement of the discrete wavelet transform, but as our framework is very general, it is easy to replace it with a transform of dictionary that suits the data best. The discrete wavelet transform is in essence a one-dimensional transform and the obvious extension to a multi-dimensional setting is not satisfactory, since it lacks directionality. The dual-tree complex wavelet transform \cite{Kingsbury1998} does have directionally selective filters and is thus better able to represent multi-dimensional directional data. We propose the following minimisation problem 
\begin{align}
&u_{\scriptscriptstyle \text{wav}} = \argmin_{u}\left\{\frac12\norm{Ku-f}_{L^2(\Sigma)}^2 +\alpha \norm{W(u)}_{L^1(\mathcal{W})} \right\} \label{eq:min_wav},
\end{align} 
where $W$ is the complex dual-tree wavelet transform and $\mathcal{W}$ is the corresponding wavelet domain. Note the close similarity in structure between \eqref{eq:min_TV} and \eqref{eq:min_wav}. Because of this similarity, we will not elaborate on the wavelet case in the numerical treatment of our models, but only mention the small changes that have to be made to the TV implementation to obtain the wavelet implementation.

\section{Numerical Framework}\label{ch:4}
In section \ref{sec:var_meth} we considered the convex, possibly non-smooth minimisation problem \eqref{eq:min_prob}. TV and TGV are defined via a maximisation problem over their dual variable(s), which can be seen in \eqref{eq:def_TV} and \eqref{eq:def_TGV}. Substituting the definition of TV/TGV in the minimisation problem changes it into a saddle-point problem, which can be solved through the use of primal-dual algorithms. These primal-dual algorithms can also be used for our wavelet regulariser or other regularisers that are not defined via their dual variables: with the use of Fenchel duality, the specific problem may be rewritten to its dual maximisation problem or its saddle-point problem.

\subsection{Description of the saddle-point problem}\label{sec:SP_formulation}
With the use of Fenchel duality, we can rewrite the general minimisation problem (cf. \eqref{eq:min_prob})
\begin{equation}\label{eq:min_general}
\min_{x\in X}\left\{F(Ax)+G(x)\right\}
\end{equation}
into the general saddle-point problem
\begin{equation}\label{eq:SP_general}
\min_{x\in X}\max_{y\in Y}\left\{\langle Ax,y\rangle+G(x)-F^*(y)\right\},
\end{equation}
where $F^*$ is the convex conjugate of $F$. We obtain our TV, TGV and wavelet models by describing the functionals $F$ and $G$ and the operator $A$ appropriately, which is done in \eqref{eq:PD_TV} and \eqref{eq:PD_TGV}. 
\begin{align}\label{eq:PD_TV}
&\textbf{TV model:} ~~ y = (q,r),~~x = u,~~Ax = (Ku,\nabla u),\nonumber\\
&F_1(\tilde{q}) = \frac{1}{2}\norm{\tilde{q}-f}_{L^2(\Sigma)},~~F_2(\tilde{r}) = \alpha\norm{\tilde{r}}_{(L^1(\Omega))^d},~~G(u) = 0,\\
~~&\min_{u\in X}\max_{(q,r)\in Y}\left\{\langle Ku,q\rangle+\langle \nabla u,r\rangle+G(u)-F_1^*(q)-F_2^*(r)\right\}\nonumber.\\ &\textbf{Note: }\text{the algorithm for wavelet reconstruction is obtained by changing }\nabla \text{ to } W,\nonumber\\ &~~~~~~~~~~\text{causing }F_2\text{ to change to }F_2(\tilde{r})=\alpha\norm{\tilde{r}}_{L^1(\mathcal{W})}\nonumber\\
~\nonumber\\
&\textbf{TGV model:} ~~ y = (q,r,s),~~x = (u,v),~~Ax = (Ku,\nabla u-v,\mathcal{E}v),\label{eq:PD_TGV}\\
&F_1(\tilde{q}) = \frac{1}{2}\norm{\tilde{q}-f}_{L^2(\Sigma)},~~F_2(\tilde{r}) = \alpha\norm{\tilde{r}}_{(L^1(\Omega))^d},~~F_3(\tilde{s}) = \alpha\beta\norm{\tilde{s}}_{(L^1(\Omega))^{d\times d}},~~G_1(u) = 0,~~G_2(v) = 0,\nonumber\\
~~&\min_{(u,v)\in X}\max_{(q,r,s)\in Y}\big\{\langle Ku,q\rangle+\langle \nabla u-v,r\rangle+\langle \mathcal{E}v,s\rangle+G_1(u)+G_2(v)-F_1^*(q)-F_2^*(r)-F_3^*(s)\big\}\nonumber.
\end{align}
Here we made a splitting $F(\cdot) = \sum_i F_i(\cdot)$ and $G(\cdot) = \sum_i G_i(\cdot)$, in the same manner as in \cite{Sidky2012}. Note that we use tildes on the operators that appear in their convex conjugate form in the primal-dual problem. The convex conjugate, also known as Fenchel–Legendre transform, is defined as
\begin{equation}\label{eq:conj_F}
F^*(x) = \max_{\hat{x}}\left\{\langle x,\tilde{x}\rangle-F(\tilde{x})\right\}.
\end{equation}

\subsection{Proximal operators for TV, TGV and wavelet model}\label{sec:Prox_ops}
A key tool for the use of primal-dual algorithms is the so called proximal splitting method. The proximal operator for a scaled functional $\gamma F(x)$ is defined as
\begin{eqnarray}\label{eq:prox}
\prox_{\gamma F}(x) = \argmin_y\left\{\frac12\norm{x-y}^2_2+\gamma F(y)\right\},
\end{eqnarray}
which gives a small value of $\gamma F(y)$, while keeping $y$ close to $x$. Because the right hand side consists of a strongly convex quadratic part and a convex functional $\gamma F(y)$, it is strongly convex and thus \eqref{eq:prox} has a unique solution.

As can be seen in section \ref{sec:PD_algorithm}, the saddle-point problems \eqref{eq:PD_TV} and \eqref{eq:PD_TGV} now have the correct structure to be numerically implemented. Therefore, the proximal operators for all $F^*_i$ and $G_i$ are needed, which are derived below.  
The convex conjugate for $F_1$ is 
\begin{equation}\label{eq:conj_F1}
F_1^*(q) = \max_{\hat{q}}\left\{\langle q,\tilde{q}\rangle-\frac12\norm{\tilde{q}-f}_{L^2(\Omega)}^2\right\}.
\end{equation}
After checking the first and second order optimality conditions, we find that $\hat{q} = f+\tilde{q}$ maximises \eqref{eq:conj_F1}. We fill this in and calculate the proximal operator
\begin{align}
\prox_{\gamma F_1^*}(q) &= \argmin_{\tilde{q}}\left\{\frac12\norm{\tilde{q}-q}_{L^2(\Omega)}^2+\gamma F_1^*(\tilde{q})\right\} = \frac{q-\gamma f}{1+\gamma}.
\end{align}
The convex conjugate for $F_2$ is
\begin{align}\label{eq:conj_F2}
F_2^*(r) &= \max_{\tilde{r}}\left\{\langle r,\tilde{r}\rangle-\alpha\norm{\tilde{r}}_{(L^1(\Omega))^d}\right\}=\begin{cases}0 &\text{ if } |\tilde{r}|_2\leq \alpha \text{ pointwise},\\
\infty &\text{ if } |\tilde{r}|_2>\alpha \text{ pointwise.}\end{cases}
\end{align}
Here we choose the pointwise 2-norm of $r$, which means $|r|_2 = |(r_1,r_2,\dots,r_d)|_2 = \sqrt{r_1^2+r_2^2+\dots+ r_d^2}$ for $r\in (L^1(\Omega))^d$. This choice has the effect that the TV-eigenfunction will be a $d$-dimensional sphere instead of a $d$-dimensional cube (1-norm) or a $d$-dimensional diamond ($\infty$-norm) \cite[Theorem 4.1]{Esedoglu2004}. For the wavelet model, $\tilde{r}$ is a scalar function instead of a vector function, which means that all pointwise $p$-norms are equal. The proximal operator for $F^*_2$ is calculated
\begin{align}
\prox_{\gamma F_2^*}(r) &= \argmin_{\tilde{r}}\left\{\frac12\norm{\tilde{r}-r}_{L^2(\Omega)}^2+\gamma F_2^*(\tilde{r})\right\} = \frac{\alpha r}{\max\{\alpha,|r|_2\}}.
\end{align}
Similarly, the proximal operator for $F^*_3$ reads
\begin{equation}
\prox_{\gamma F_3^*}(s) = \frac{\alpha\beta s}{\max\{\alpha\beta,|s|_2\}}. 
\end{equation}
It is easily seen that $\prox_{\gamma G_1}(u) = u$ and $\prox_{\gamma G_2}(v) = v$.

\subsection{First-order primal dual algorithm}\label{sec:PD_algorithm}
We make use of the modified primal-dual hybrid gradient algorithm (PDHGM) as proposed by Chambolle and Pock \cite{Chambolle2011}, which can be seen as a generalisation to the PDHG algorithm \cite{Zhu2008}. PDHGM has the advantage that it has a direct solution at every step. Furthermore it can be shown that PDHGM has a convergence rate of at least $\bigo(n)$ and can even be $O(n^2)$, where $n$ is the number of steps in the algorithm \cite{Chambolle2011, Pock2011}. 

In section \ref{sec:SP_formulation} we have split the functionals $F$ and $G$ in multiple parts. Therefore, we also need to split the proximal operators and evaluate them separately. A similar algorithm to \cite[Algorithm 4]{Sidky2012} is used.
\begin{algorithm}
\caption{First-order primal-dual algorithm for TV reconstruction}\label{alg:TV}
\begin{algorithmic}
\State \textbf{Parameters: }choose $\tau,\sigma_i>0$ s.t. $\tau\sigma_iL_i^2<1$; ~choose $\theta\in[0,1]$. 
\State \textbf{Initialisation: }set $u^0=0$, $q^0=0$, $r^0=0$, $\bar{u}^0 = 0$.\vspace{2mm}
\For{$k \gets 1$ to $N$}
\State $q^{n+1} = \prox_{\sigma_1 F_1^*}(q^n+\sigma_1K\bar{u}^n),$
\State $r^{n+1} = \prox_{\sigma_2 F_2^*}(r^n+\sigma_2\nabla\bar{u}^n),$
\State $u^{n+1} = \prox_{\tau G}(u^{n}-\tau(K^* q^{n+1}-\div ~r^{n+1}),$
\State $\bar{u}^{n+1} = u^{n+1} +\theta(u^{n+1}-u^{n}).$
\EndFor
\vspace{2mm}
\State \textbf{return} $u^N$
\end{algorithmic}
The algorithm for wavelet reconstruction is obtained by changing $\nabla$ to $W$ and $\div$ to $W^{-1}$.
\end{algorithm}

Algorithm \ref{alg:TV} shows the implementation of the TV model. We perform two separate dual steps, which are both used as input for the primal step after that. The implementation of the wavelet model is obtained by changing $\nabla$ and $\div$ to $W$ and $W^{-1}$ respectively. 

\begin{algorithm}
\caption{First-order primal-dual algorithm for TGV reconstruction}\label{alg:TGV}
\begin{algorithmic}
\State \textbf{Parameters: }choose $\tau,\sigma_i>0$ s.t. $\tau\sigma_iL_i^2<1$; ~choose $\theta\in[0,1]$. 
\State \textbf{Initialisation: }set $u^0=0$, $v^0=0$, $q^0=0$, $r^0=0$, $s^0=0$, $\bar{u}^0 = 0$, $\bar{v}^0 = 0$.
\vspace{2mm}
\For{$k \gets 1$ to $N$}
\State $q^{n+1} = \prox_{\sigma_1 F_1^*}(q^n+\sigma_1K\bar{u}^n),$
\State $r^{n+1} = \prox_{\sigma_2 F_2^*}(r^n+\sigma_2(\nabla\bar{u}^n-\bar{v}^n)),$
\State $s^{n+1} = \prox_{\sigma_2 F_3^*}(s^n+\sigma_2\mathcal{E}\bar{v}^n),$
\vspace{2mm}
\State $u^{n+1} = \prox_{\tau G_1}(u^{n}-\tau(K^* q^{n+1}-\div ~r^{n+1}),$
\State $v^{n+1} = \prox_{\tau G_2}(v^{n}-\tau(\mathcal{E}^* s^{n+1}-r^{n+1}),$
\vspace{2mm}
\State $\bar{u}^{n+1} = u^{n+1} +\theta(u^{n+1}-u^{n}),$
\State $\bar{v}^{n+1} = v^{n+1} +\theta(v^{n+1}-v^{n}),$

\EndFor
\vspace{2mm}
\State \textbf{return} $u^N$
\end{algorithmic}
\end{algorithm}

Algorithm \ref{alg:TGV} shows the implementation of the TGV model, which is built in the same spirit as the TV algorithm. 

\subsubsection{Discretisation}
The measurement data $f$ is sampled in both space and time: the detectors are modelled as point detectors at specific locations and measure the pressure at sampled time instances with interval $\Delta t$. The image that we wish to reconstruct is also a discretised one. In case these discretisation do not perfectly match (which is often the case), we might get significant errors. To avoid these errors, we use the interpolation procedure as explained in \cite[Chapter 5]{Willemink2010}.

Since we are using an iterative procedure, it is computationally expensive to evaluate the discrete integral operator $K$ at every step. Therefore, we calculate the discretised matrix $K$ once and use it iteratively. Instead of analytically determining $K^*$ and then discretising it, we use the adjoint matrix $K^\textsc{T}$. Accordingly, we also use matrix versions for the operators $\nabla$ and $\mathcal{E}$. For the complex dual-tree wavelet transform and its adjoint (which in this case is equivalent to its inverse), we make use of the implementation by Cai and Li\footnote{available at \url{http://eeweb.poly.edu/iselesni/WaveletSoftware/dt2D.html}}.

\subsubsection{Algorithm parameter selection}
In the work of Chambolle and Pock \cite{Chambolle2011}, it is shown that the PDHGM algorithm always converges if $\sigma\tau\norm{A}^2<1$. Here $A$ is the combined operator as defined in \eqref{eq:PD_TV} and \eqref{eq:PD_TGV}. A uniform bound for $\tau$ and $\sigma$ might be undesired if the separate operators $K$ and $\grad$, $\mathcal{E}$ or $W$ have very different norms: all the proximal steps have a small step size, while this might not be needed for the steps related to only one of these operators.

In \cite{Pock2011} this problem has been solved for matrix operators. Instead of using one value for $\sigma$ and $\tau$ in the PDGHM algorithm, one can use diagonal matrices $\Sigma$ and $T$ that satisfy $\norm{\Sigma^{1/2}AT^{1/2}}<1$. With this choice, the sequence generated by the algorithm weakly converges to an optimal solution of the saddle-point problem. 

For our specific operators, it appears that $\norm{K}\gg\norm{\nabla}\approx\norm{\mathcal{E}}\approx\norm{W}=1$, but they show a very similar structure within each operator. Therefore, we only choose two sets of two operators, since they have shown to give a smoother convergence plot than with the parameter choice in \cite{Pock2011}, while the convergence rate is similar. We choose one set of two parameters for the operator $K$ and one for the combination of other operators. More precisely, we choose $\sigma_1,\sigma_2,\tau$ such that $\sigma_1\tau L_1 ^2<\frac14$ and $\sigma_2\tau L_2 ^2<\frac14$. Here $L_1 = \norm{K}$, $L_2 = \norm{\nabla}$ or $L_2 = \norm{W}$ in the TV or wavelet model and $L_2$ is the norm of the combined other operators in the TGV model.

It is easily shown that the bound $\norm{\Sigma^{1/2}AT^{1/2}}<1$ is satisfied by this choice if we write $A$ as a concatenation of the matrix version of all operators. For the TGV model, we have
\begingroup
\arraycolsep=2.5pt
\begin{equation*}
A = \left[
\begin{array}{c c}
K & 0 \\ \hdashline[2pt/3pt]
\nabla & -I \\
0 & \mathcal{E}
\end{array}
\right]\in\R^{(m_1+m_2)\times n},
\end{equation*}
\endgroup
and the diagonal matrices
\begin{align*}
\Sigma = \begin{bmatrix}
\Sigma_1 & 0\\
0 & \Sigma_2 \end{bmatrix}
= \begin{bmatrix}
\sigma_1 I_{m_1} & 0\\
0 & \sigma_2 I_{m_2} \end{bmatrix},~~T = \tau I_n.
\end{align*}
This gives us the estimate
\begin{align}
\begin{aligned}
\norm{\Sigma^{1/2}AT^{1/2}} \leq &\norm{\Sigma_1^{1/2}
\begin{bmatrix}
K & 0 
\end{bmatrix} T^{1/2}} + \norm{\Sigma_2^{1/2}
\begin{bmatrix}
\nabla & -I \\
0 & \mathcal{E}
\end{bmatrix} T^{1/2}}\\
= &\sqrt{\sigma_1\tau} L_1 + \sqrt{\sigma_2\tau} L_2<\frac12+\frac12 = 1.\end{aligned}
\end{align}
A similar estimate can be given for the TV and wavelet model.

\subsubsection{Algorithm validation}
The PDHGM algorithm with its specific discretisation and parameter choices is validated by looking at two criteria. Firstly, the duality gap is taken into consideration. The duality gap gives the (non-negative) difference between the primal and the dual functional. As the duality gap approaches zero, the solution gets asymptotically close to the desired solution for the primal-dual problem \cite{Boyd2004}. Hence the algorithm is validated by checking if the duality gap approaches zero. Moreover, a stopping criterion for the algorithm can be constructed based on the duality gap. In the formulation of \eqref{eq:SP_general}, the duality gap reads
\begin{equation}
D^n = F(Ax^n)+G(x^n)+F^*(y^n)+G^*(-A^*y^n).
\end{equation}
The second criterion is the residual of the variables, defined as
\begin{equation}
x^n_\text{res} = \frac{\norm{x^n-x^{n-1}}}{x^{n-1}},
\end{equation}
where any variable can be filled in instead of $x$. The residual $x^n_\text{res}\to0$ for $n\to\infty$ whenever $(x^n-x^{n-1})\to0$, which means that the sequence $x^n$ converges. 

\section{Experimental setup}\label{sec:exp_setup}
The proposed methods are defined for solving the general problem $Ku=f$. Here $u$ is the desired reconstruction, $f$ is (preprocessed) data and $K$ is some operator that maps $u$ to $f$. Because of this generic structure, the methods can be applied to reconstruct images for many photoacoustic tomographs with many different forward models.

\begin{figure*}[!ht]
\centering
\includegraphics[width=0.9\linewidth]{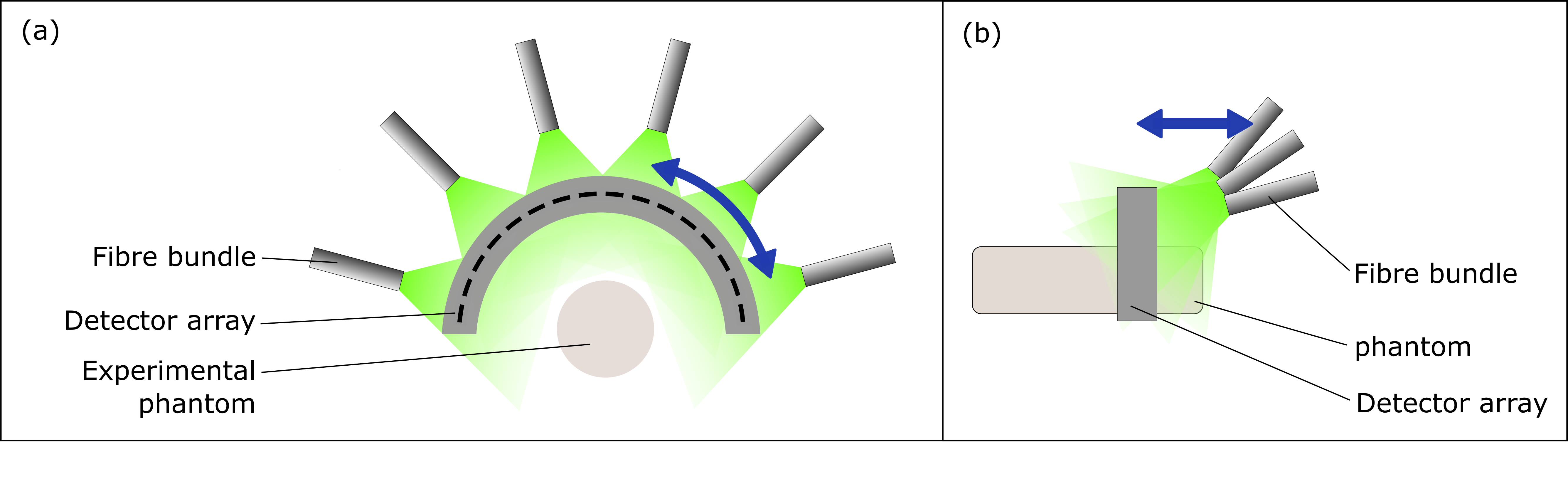}
\caption{Schematic overview of the experimental setup, where a phantom is measured.}
\label{fig:exp_setup}
\end{figure*}

In our experimental setup, we make use of the finger joint imager as specified in \cite{Es2015}. A side-illumination laser delivers pulses of optical energy with a wavelength of 532 nm. For the recording of pressure waves a 1D piezoelectric detector array with 64 elements in a half-circle is used. The detector array has a central frequency of 7.5 MHz. Both the detectors and fibre bundles are simultaneously rotatable over 360 degrees. They can also be translated in order to image multiple slices. The detectors have a narrow focus (0.6 mm plane thickness) in one dimension, making it suitable for 2D slice based imaging. A schematic overview of the experimental setup is shown in Figure \ref{fig:exp_setup}. For a more extensive explanation of the setup, we refer to \cite{Es2015}.

\subsection{Forward model}
Based on the work of \cite{Kruger1995} and \cite{Wang2004}, we write our forward model as the projection of the absorbed energy distribution $u(\mathbf{x})$ over spherical surfaces. Under the assumption of a homogeneous speed of sound $c_0$, we obtain the following expression for the pressure measured by the detector at location $\mathbf{x}$:
\begin{eqnarray}
p(\mathbf{x},t) = \frac{\beta}{4\pi C_p}\left(\frac{1}{t}\iint_{|\mathbf{x}-\tilde{\mathbf{x}}|=c_0t} u(\tilde{\mathbf{x}}) \text{d} \tilde{\mathbf{x}}\right)\ast_t\frac{\partial I(t)}{\partial t}\ast_t h_{IR}(t),
\end{eqnarray}
where $\beta$ is the thermal expansion coefficient, $C_p$ is the specific heat, $I(t)$ is the illumination profile of the laser, $h_{IR}(t)$ is the impulse response of the detectors and $u(\mathbf{x})$ is the absorbed energy. Instead of finding an explicit expression for $\frac{\partial I(t)}{\partial t}$ and $h_{IR}(t)$, we make use of the calibration measurement $p_\text{cal}$ of an approximate delta pulse \cite{Wang2004} located at $\mathbf{x}_p$, i.e. $u_\text{cal}(\mathbf{x}) \approx \delta(\mathbf{x}-\mathbf{x}_p)$. After finding an analytic expression for $p_\text{cal}(\mathbf{x},t)$ in terms of $\frac{\partial I(t)}{\partial t}$ and $h_{IR}(t)$ \cite{Willemink2010}, we obtain
\begin{equation}\label{eq:FW_full}
p(\mathbf{x},t) = |\mathbf{x}-\mathbf{x}_p|\left(\frac{1}{t}\iint_{|\mathbf{x}-\tilde{\mathbf{x}}|=c_0t} u(\tilde{\mathbf{x}})\text{d} \tilde{\mathbf{x}}\right)\ast_t p_\text{cal}\left(\mathbf{x},t'\right),
\end{equation}
where $t'$ is the retarded time $t-\frac{|\mathbf{x}-\mathbf{x}_p|}{c}$.  

\subsection{Preprocessing for reconstruction}\label{sec:prep}
We can write \eqref{eq:FW_full} concisely as 
\begin{align}
&p = \left(\frac{|\mathbf{x}-\mathbf{x}_p|}{t}\left( Ku\right)\right)\ast p_\text{cal},\nonumber\\
\text{with } Ku &:= \iint_{|\mathbf{x}-\tilde{\mathbf{x}}|=c_0t} u(\tilde{\mathbf{x}})\diff\tilde{\mathbf{x}}.\label{eq:K_def}
\end{align}
We left the dependencies on $\mathbf{x}$, $t$ and $t'$ out for an uncluttered notation; the convolution with $\tilde{p}_\text{cal}$ is performed in the time domain. Recall that we are solving the problem $f = Ku$, where in this specific setup, $K$ is defined as in \eqref{eq:K_def}. In order to get our preprocessed data $f$, we first solve the deconvolution problem
\begin{align}
p = f_t\ast p_\text{cal}, \label{eq:FW_conv}\\
\text{with }f_t := \frac{|\mathbf{x}-\mathbf{x}_p|}{t}f.\label{eq:f_rewrite}
\end{align}
After applying the convolution theorem to \eqref{eq:FW_conv}, we obtain $\F\{p\} = \F\{f_t\}\F\{p_\text{cal}\}$, which for $\F\{p_\text{cal}\}>0$ is equivalent to 
\begin{equation}\label{eq:conv_thm2}
f_t = \F^{-1}\left\{\frac{\F\{p\}}{\F\{p_\text{cal}\}}\right\}.
\end{equation}
Unfortunately, the right hand side of \eqref{eq:conv_thm2} is undefined when $\F\{p_\text{cal}\}=0$. Moreover any noise on $p_\text{cal}$ can have a big influence on the expression. Therefore, we construct a regularised deconvolution filter
\begin{equation*}
h_\text{dec} := \F^{-1}\left\{\frac{\overline{\F\{p_\text{cal}\}}}{|\F\{p_\text{cal}\}|^2+\varepsilon}\right\},
\end{equation*}
such that 
\begin{equation}\label{eq:deconv}
p\ast h_\text{dec} = \F^{-1}\left\{\frac{\F\{p\}\overline{\F\{p_\text{cal}\}}}{|\F\{p_\text{cal}\}|^2+\varepsilon}\right\} \approx f_t.
\end{equation}
Here $\varepsilon$ is a small parameter such that it suppresses noise on $p_\text{cal}$, but leaves higher parts of it intact. By consecutively applying \eqref{eq:deconv} and \eqref{eq:f_rewrite}, we obtain our preprocessed measurement data for \eqref{eq:objective}.

\section{Test objects}\label{sec:phantoms}
In order to test our reconstruction framework with its regularisers, both digital and experimental phantoms have been created. For any of the three regularisers a specific digital phantom has been created, that fulfils the prior assumption on the image to be reconstructed. An experimental phantom with a vascular structure has been created to test the effectiveness of the different regularisers on real data.

\subsection{Digital phantoms}
\begin{figure}[!ht]
\centering
\begin{subfigure}{0.3\linewidth}
\includegraphics[width=\linewidth]{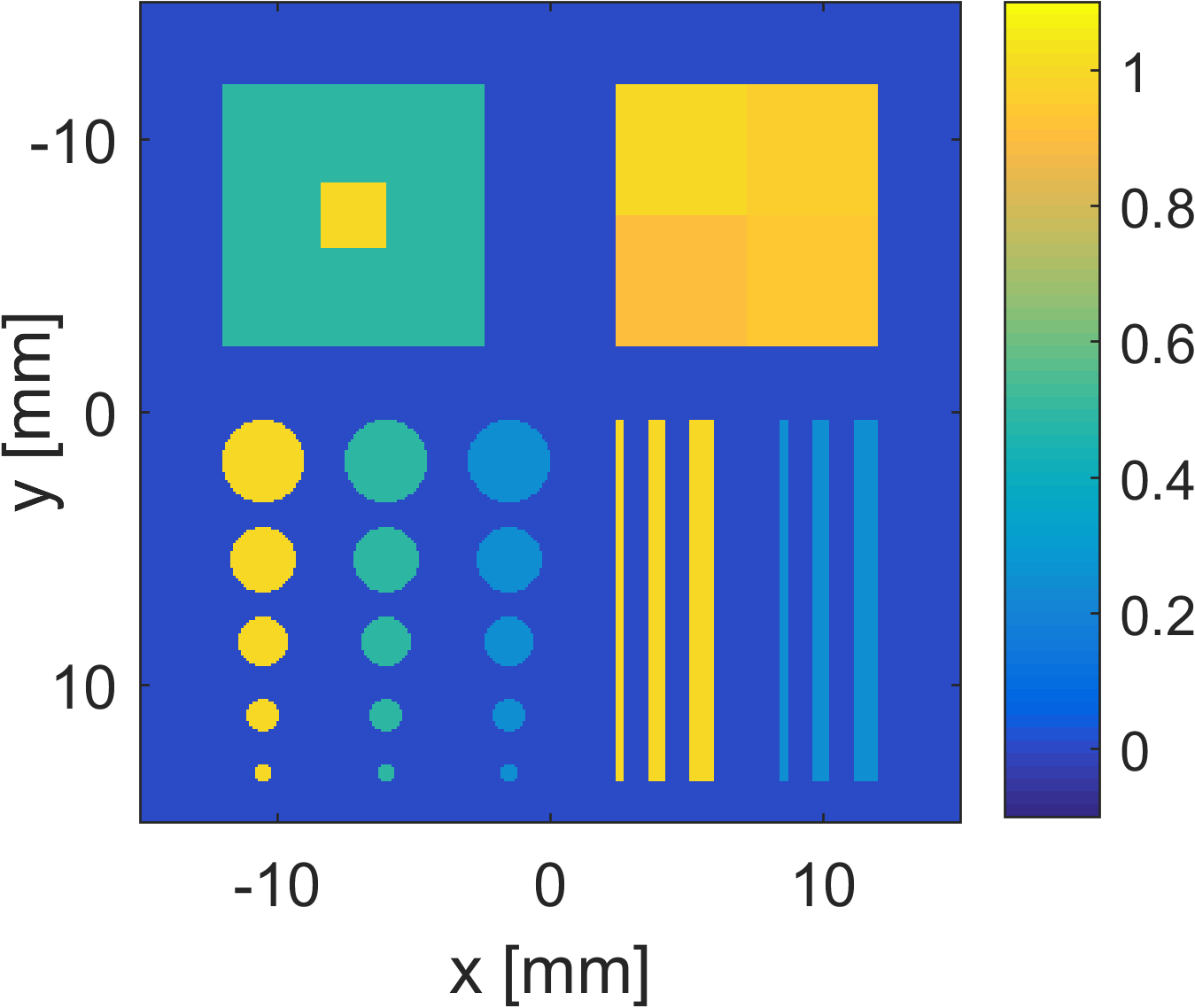}
\caption{}
\label{fig:synth_diverse}
\end{subfigure}
\begin{subfigure}{0.30\linewidth}
\includegraphics[width=\linewidth]{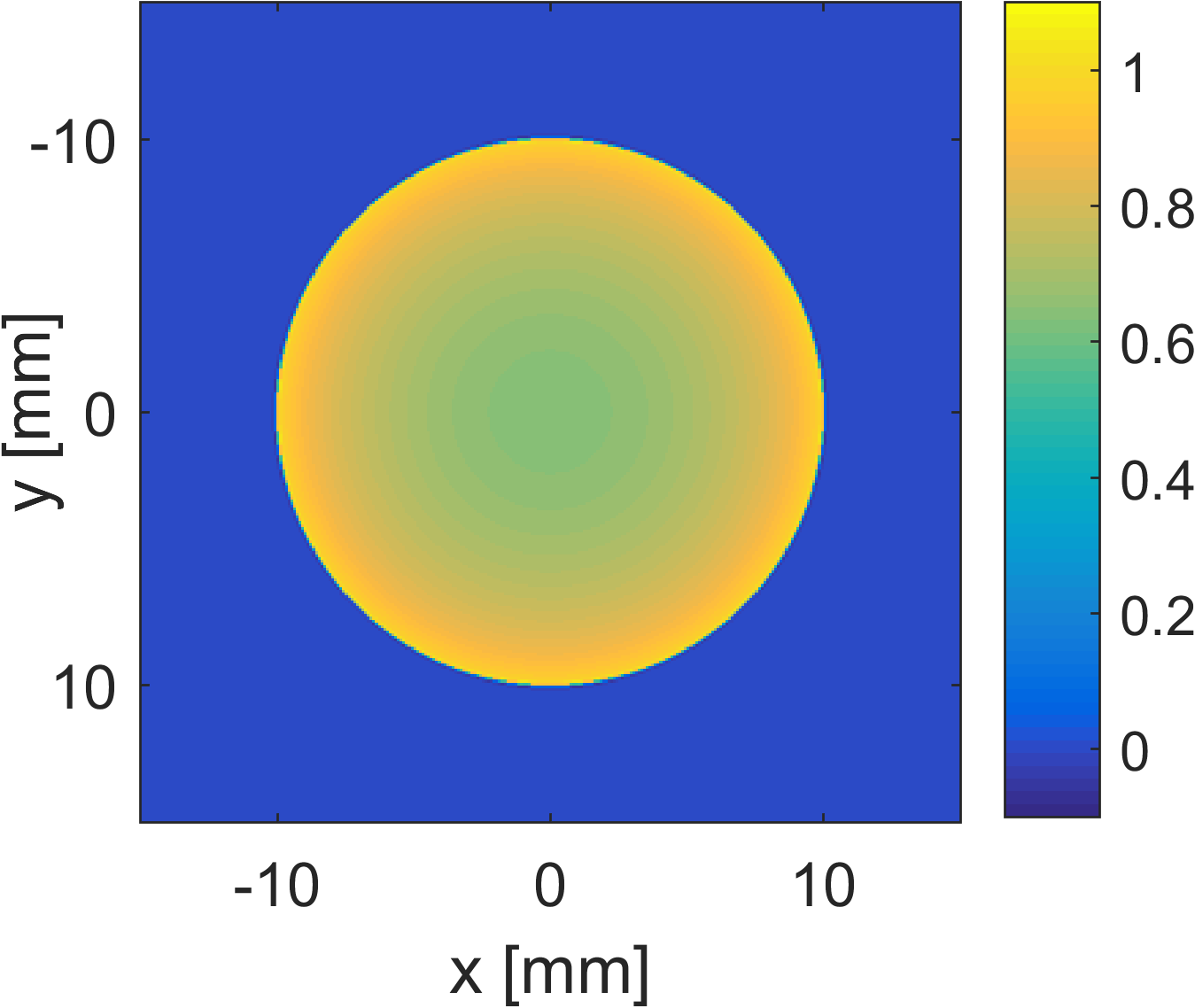}
\caption{}
\label{fig:synth_exp}
\end{subfigure}
\begin{subfigure}{0.30\linewidth}
\includegraphics[width=\linewidth]{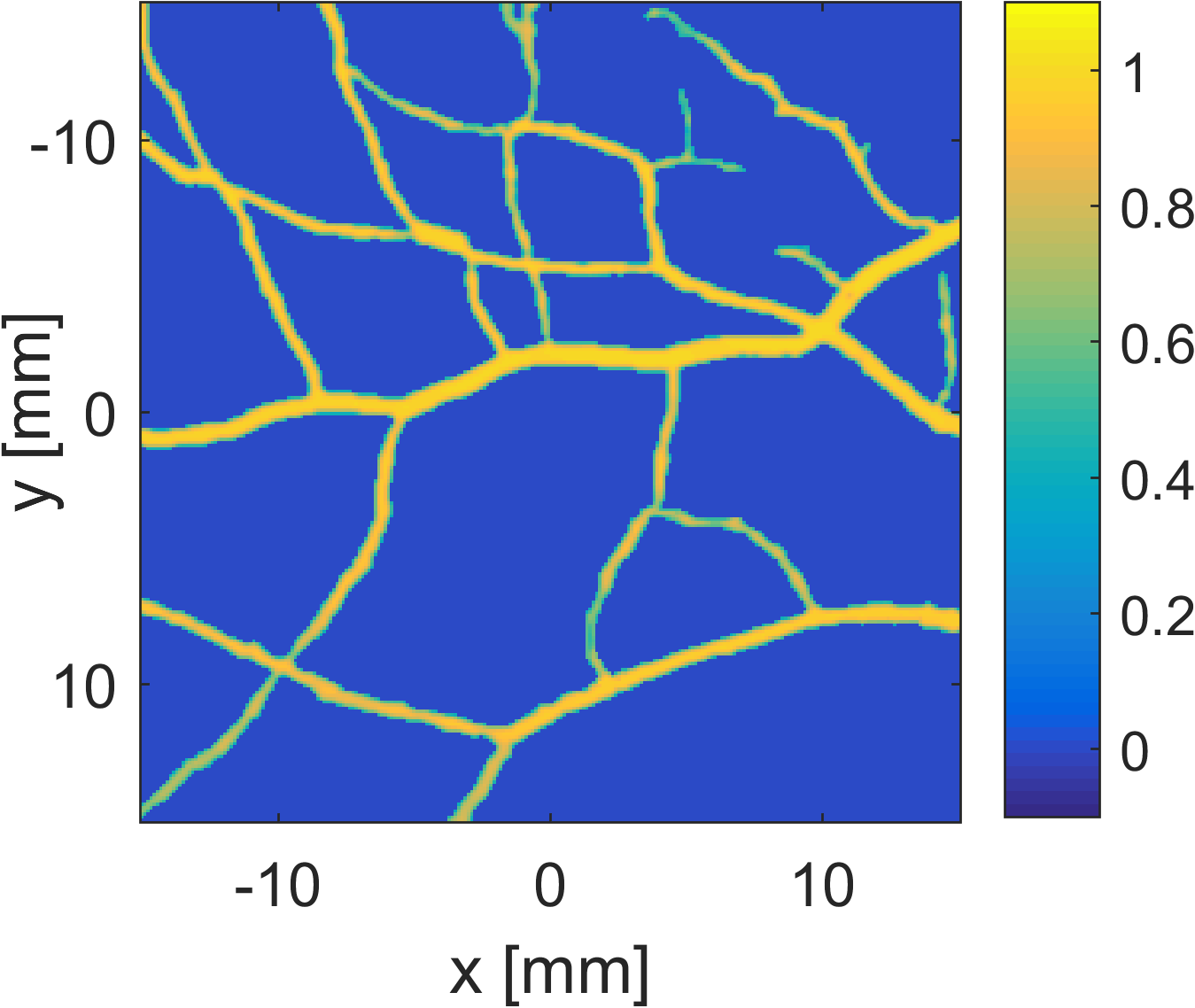}
\caption{}
\label{fig:synth_vascular}
\end{subfigure}
\caption{(a) Ground truth digital phantom containing multiple piecewise constant structures of various shapes, sizes and intensities. (b) Digital phantom containing disc with smoothly decaying intensity from outside to inside. (c) Digital phantom containing a vascular structure.}
\label{fig:synth}
\end{figure}

The first digital phantom is shown in Figure \ref{fig:synth_diverse}. For this phantom, it is assumed that the fluence rate in the illuminated object is homogeneous. This means that sharp changes in the absorbed energy are expected and it thus consists of piecewise constant values. This image consists of both large and small objects with a variety of shapes, among which squares, discs and elongated rectangles. The larger discs could resemble dermis in the finger, while the small discs and rectangles could resemble blood vessels respectively perpendicular and in the 2D-plane of focus. Because of the sharp discontinuities, we will reconstruct this image with the TV method and compare this with FBP and LS-T reconstructions.

The second digital phantom is shown in Figure \ref{fig:synth_exp}. For this phantom, it is assumed that the fluence rate is decaying from the edge of the disc to the middle of the disc. It resembles tissue that absorbs enough light such that the absorbed energy decays as we are deeper inside the tissue. Since the fluence rate is assumed heterogeneous, the TGV method will be used for the reconstruction.

The third digital phantom, shown in Figure \ref{fig:synth_vascular}, is a vascular structure, which has been obtained by preprocessing part of a retinal image from the DRIVE database \cite{Staal2004}. Such a structure can be expected when blood vessels lie in the plane of focus. Moreover, such a phantom is similar to 3D vascular structures that can be expected in photoacoustic breast imaging. 

\subsection{Synthetic data acquisition}
The acquisition of the preprocessed data $f$ from the measured pressure $\tilde{p}$ has a strong dependency on the calibration measurement $\tilde{p}_\text{cal}$, as can be seen in section \ref{sec:prep}. If we want to simulate the acquisition of $\tilde{p}$, we will also have to make use of a calibration measurement $\tilde{p}_\text{cal}$, as can be seen in \eqref{eq:FW_full}. In order to avoid an inverse crime, we will use a different calibration measurement for the forward model as for the reconstruction. Moreover, the discretisation of the ground truth image $u_{GT}$ will be different from the discretisation of the reconstructed image $\hat{u}$.

\subsection{Experimental phantoms}
In order to experimentally test the capability of the regularisers, we developed a phantom with absorbers that resemble a vascular structure. The vessel-shaped absorbers (filaments) were constructed of sodium alginate (SA) gel carrying iron oxide nanoparticles to impart absorption. For the latter we used commercially available superparamagnetic iron oxide (SPIO) nanoparticles (Endorem - Guerbet, Villepinte, France). A dilute Endorem dispersion was mixed with SA solution in distilled water to arrive at a final 2\% (w/v) of SA solution. The filaments were fabricated by extruding the SA solution through a syringe with a 30g needle and allowing to fall into 0.7 M calcium chloride (CaCl$_2$) solution. SA undergoes gel formation in the presence of Ca$^{2+}$ ions in water. The resulting gel sinks to the bottom and hardens for 15 minutes. Finally, the filaments are isolated and washed three times with distilled water to ensure removal of residual Ca$^{2+}$.

The cylindrical phantom (diameter 24 mm) was made of Agar gel with Intralipid to provide optical scattering. To create the gel, first 3\% (v/v) Agar was dissolved in water by boiling it in the microwave until a clear solution had been formed. Next, the temperature of the mixture was decreased to 40$^{\circ}$C under continuous stirring. A 3\% (v/v) aqueous solution of 20\% stock Intralipid was added drop-wise with stirring. This provides a reduced scattering coefficient ($\mu$s') of 9.7 cm$^{-1}$ at 532 nm.

\begin{figure}[!ht]
\centering
\begin{subfigure}{0.3\linewidth}
\includegraphics[width=\linewidth]{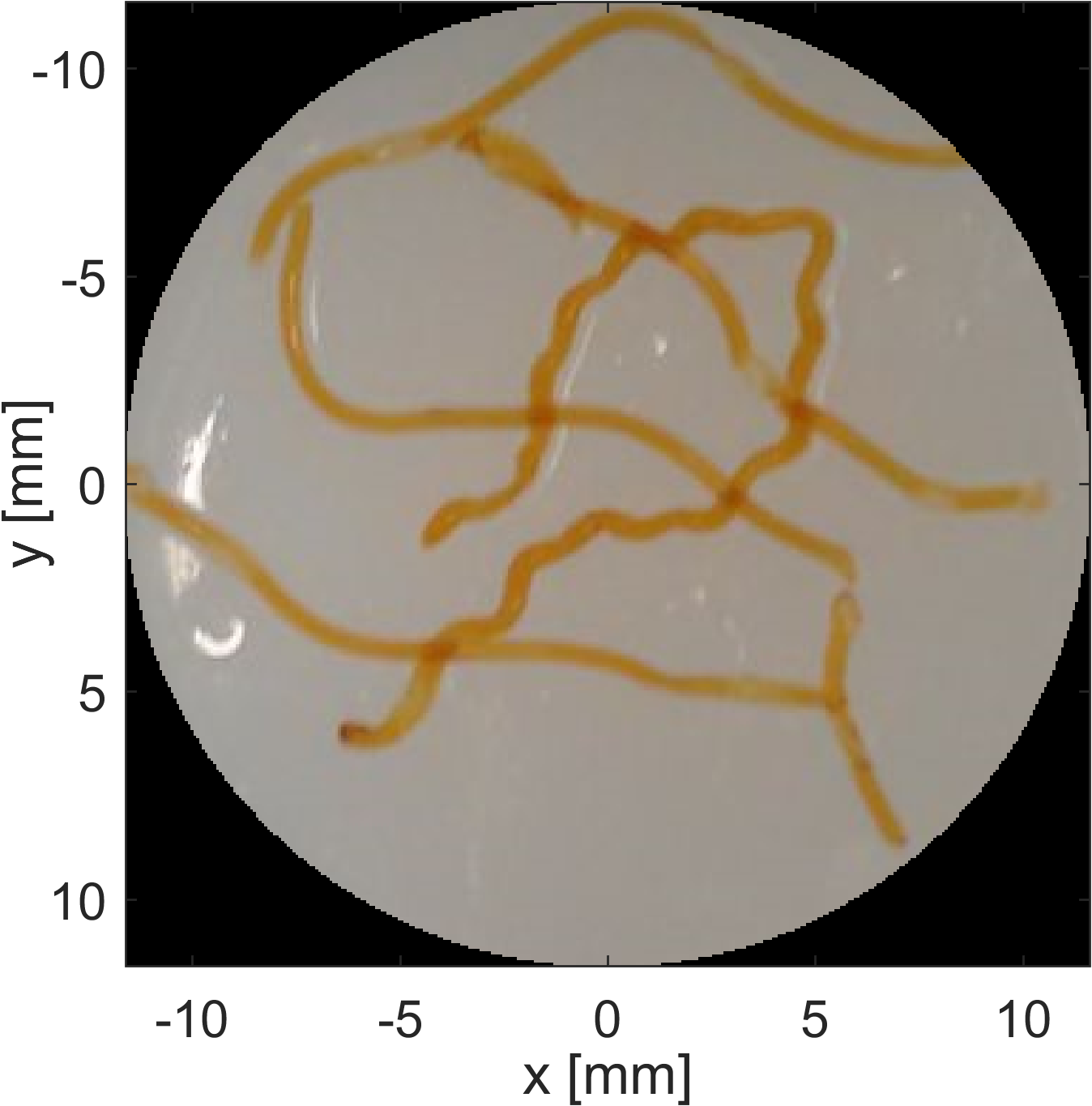}
\caption{}
\label{fig:exp_original}
\end{subfigure}
\begin{subfigure}{0.30\linewidth}
\includegraphics[width=\linewidth]{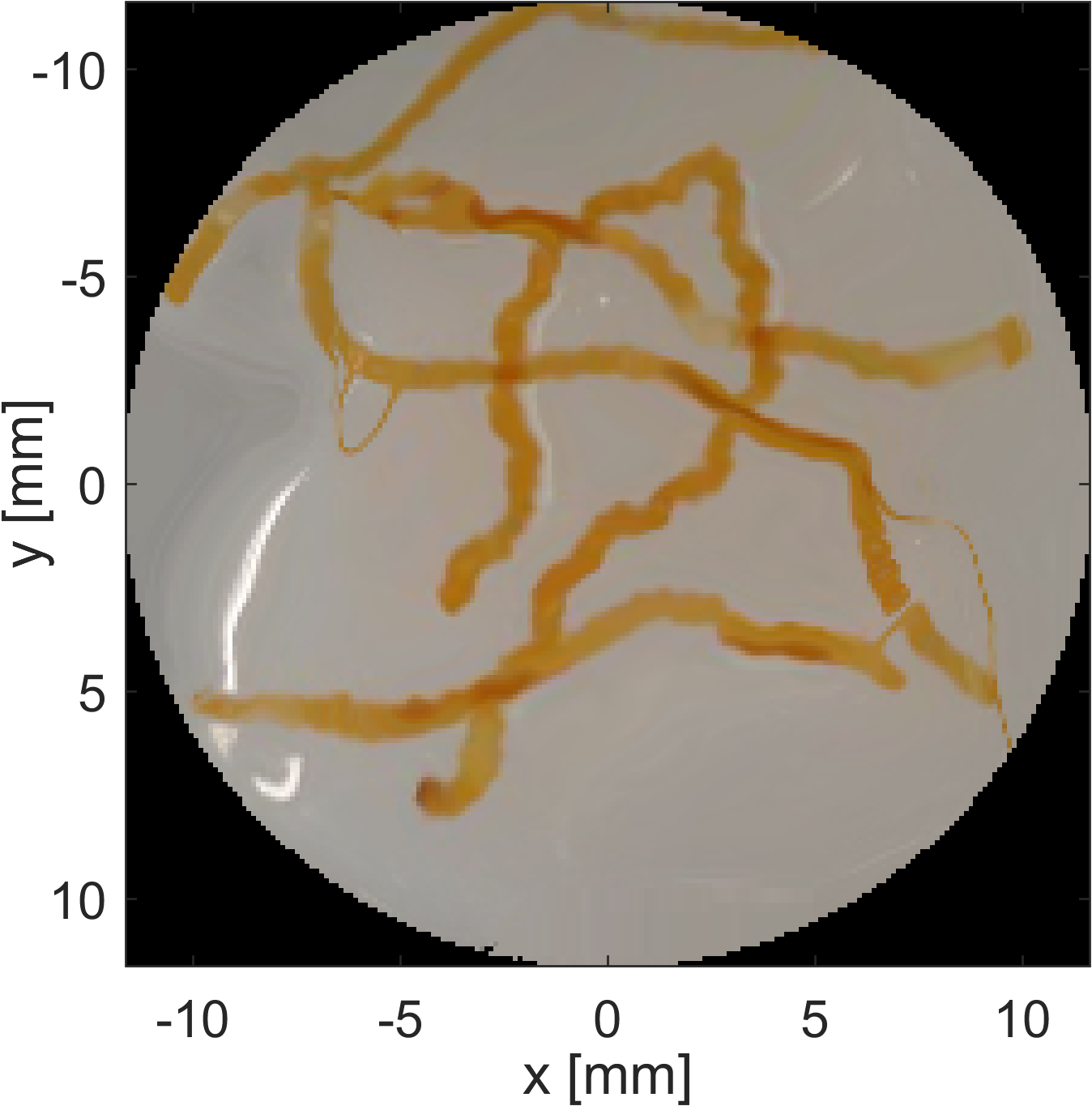}
\caption{}
\label{fig:exp_registered}
\end{subfigure}
\begin{subfigure}{0.30\linewidth}
\includegraphics[width=\linewidth]{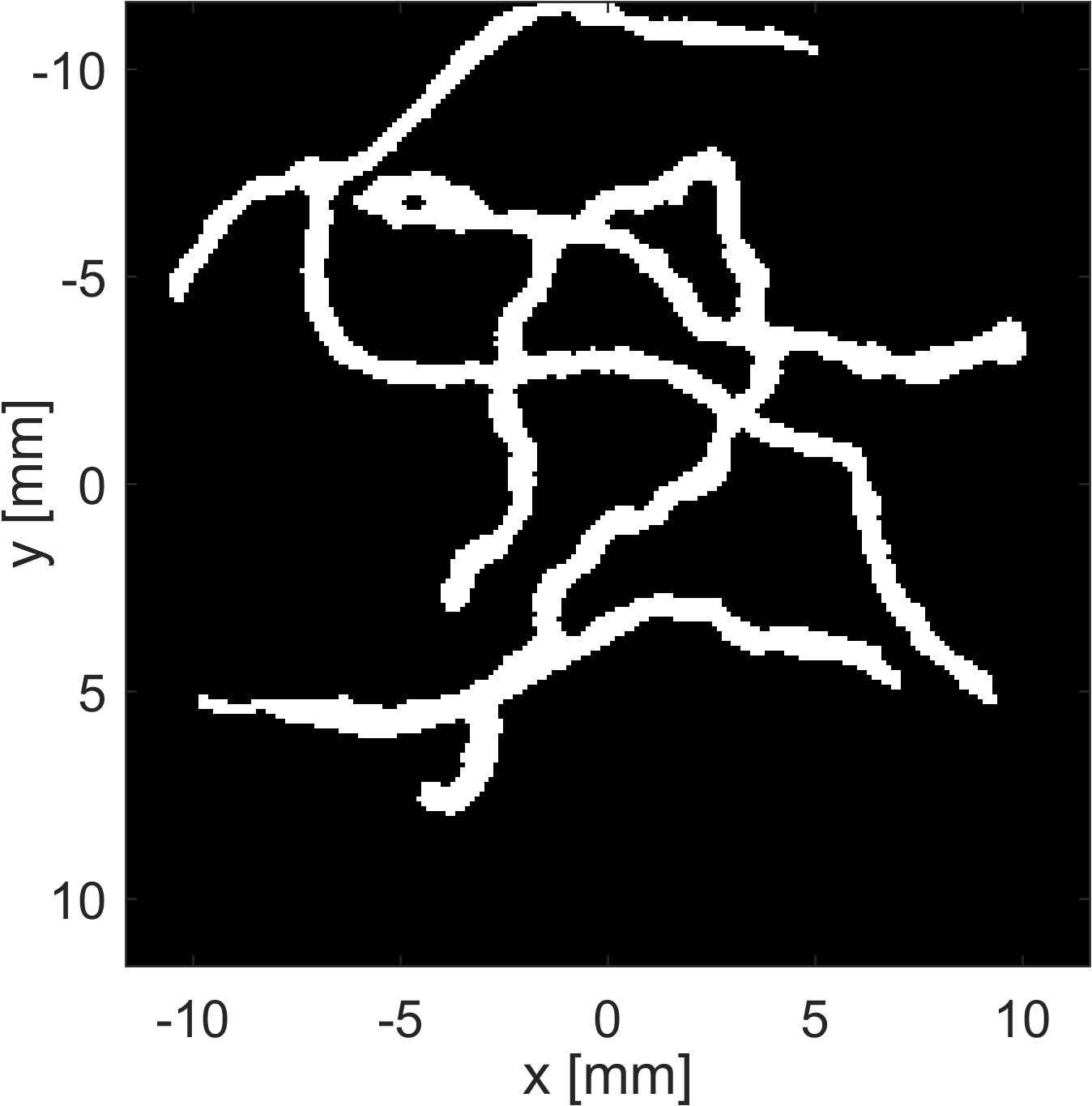}
\caption{}
\label{fig:exp_segmented}
\end{subfigure}
\caption{(a) Experimental phantom containing multiple anisotropic vascular-like structures. (b) Registered experimental phantom (aligned with photoacoustic reconstructions). (c) Thresholded registered experimental phantom.}
\label{fig:exp}
\end{figure}

To prepare the phantom, the Agar solution was poured halfway into a tube as mould and allowed to harden. A certain number of filaments were then laid on the surface of the stiff Agar gel (Figure \ref{fig:exp_original}). The whole was then topped up with more Agar solution by carefully pouring. When the absorbing agar solution had set, the phantom was removed from the mould.

As the reconstructions will later show, some of the filaments moved during the pouring of the Agar solution. For later comparisons, we have registered (aligned) the photograph of the experimental phantom with the reconstruction. This has been achieved with an affine and B-spline grid based registration, for which a toolbox is available on \textsc{Matlab} Central\footnote{Available at \href{https://nl.mathworks.com/matlabcentral/fileexchange/20057-b-spline-grid--image-and-point-based-registration}{\tt https://nl.mathworks.com/matlabcentral/fileexchange/20057-b-spline-grid{-}-image-and-point-based-registration}}. The registered image is shown in Figure \ref{fig:exp_registered}. The registered image has been segmented to enable future analysis. The segmented image is shown in Figure \ref{fig:exp_segmented}. The thin lines in Figure \ref{fig:exp_registered} are registration artefacts and have been removed. The segmentation serves as a (quasi) ground truth reference for validating the morphology of our experimental reconstructions.

\subsection{Quality measures}
For the synthetic phantoms, a digital ground truth image is available to compare the reconstructions with. Besides a visual comparison, we make use of the peak signal-to-noise ratio (PSNR), defined as 
\begin{eqnarray*}
PSNR(\hat{u},u_{GT}) = 10 \log_{10}\left(\frac{\max(u_{GT})}{\norm{\hat{u}-u_{GT}}_2^2/N}\right),
\end{eqnarray*}
where $\hat{u}$ is the reconstructed image, $u_{GT}$ the ground truth image and $N$ the number of pixels in the image.

A quality measure based on image intensities is not possible for the experimental reconstructions, since no digital ground truth image is available. However, when interested in the geometry and location of the vascular structure, only the segmentation of the reconstruction is of importance. This reconstruction segmentation can be compared with the ground truth segmentation (Figure \ref{fig:exp_segmented}). We follow this idea and make use of the receiver operating characteristic (ROC) curve. The ROC curve is obtained by plotting the false positive rate (one minus specificity) against the true positive rate (sensitivity) of the thresholded reconstruction for various threshold values. The true positive rate gives the fraction of pixels inside the ground truth segmentation that are correctly classified as such by the reconstruction segmentation. The false positive rate gives the fraction of pixels outside the ground truth segmentation that are incorrectly classified as such by the reconstruction segmentation. As the threshold for the reconstruction varies, more pixels will be correctly segmented, at the cost of incorrectly segmented pixels. An ideal situation would be one where the true positive rate increases, without changing the false positive rate, i.e. an almost vertical line, followed by an almost horizontal line.

\section{Synthetic results}\label{sec:results}
Our reconstruction framework has been tested on both synthetic and experimental data. It is compared with two standard direct reconstruction methods, namely filtered backprojection (FBP) and Tikhonov-regularised least squares (LS-T), cf. section \ref{sec:var_meth}. The specific FBP algorithm that was used is explained in \cite{Willemink2010}. In the synthetic case, for the LS-T, TV, TGV and wavelet methods, the regularisation parameters are chosen such that the PSNR between the ground truth and the reconstruction are best. In case of data with additive Gaussian noise, the optimal regularisation parameters can be found with an L-curve method \cite{Hansen2000}. Robustness against noise and compressive sampling will be shown.

\subsection{Robustness under compressive sampling}
Synthetic preprocessed measurement data have been obtained according to \eqref{eq:FW_full} with a uniform sampling of 64 dectors over a 172 degree array. Six rotations of 60 degrees have been performed, giving us an almost uniform sampling in a total of 384 detection locations. For the first comparison, different reconstruction methods under uniform compressive sampling are considered. It is unclear if uniform compressive sampling gives the best reconstruction quality with respect to the different reconstruction methods (see e.g. Haber \textit{et al.} on experimental design \cite{Haber2008}). However, it intuitively makes sense to place the detectors uniformly. Moreover, since the detector need some space, they are often placed in such a fashion. 

\begin{figure*}[!htb]
\centering
\includegraphics[width=0.9\linewidth]{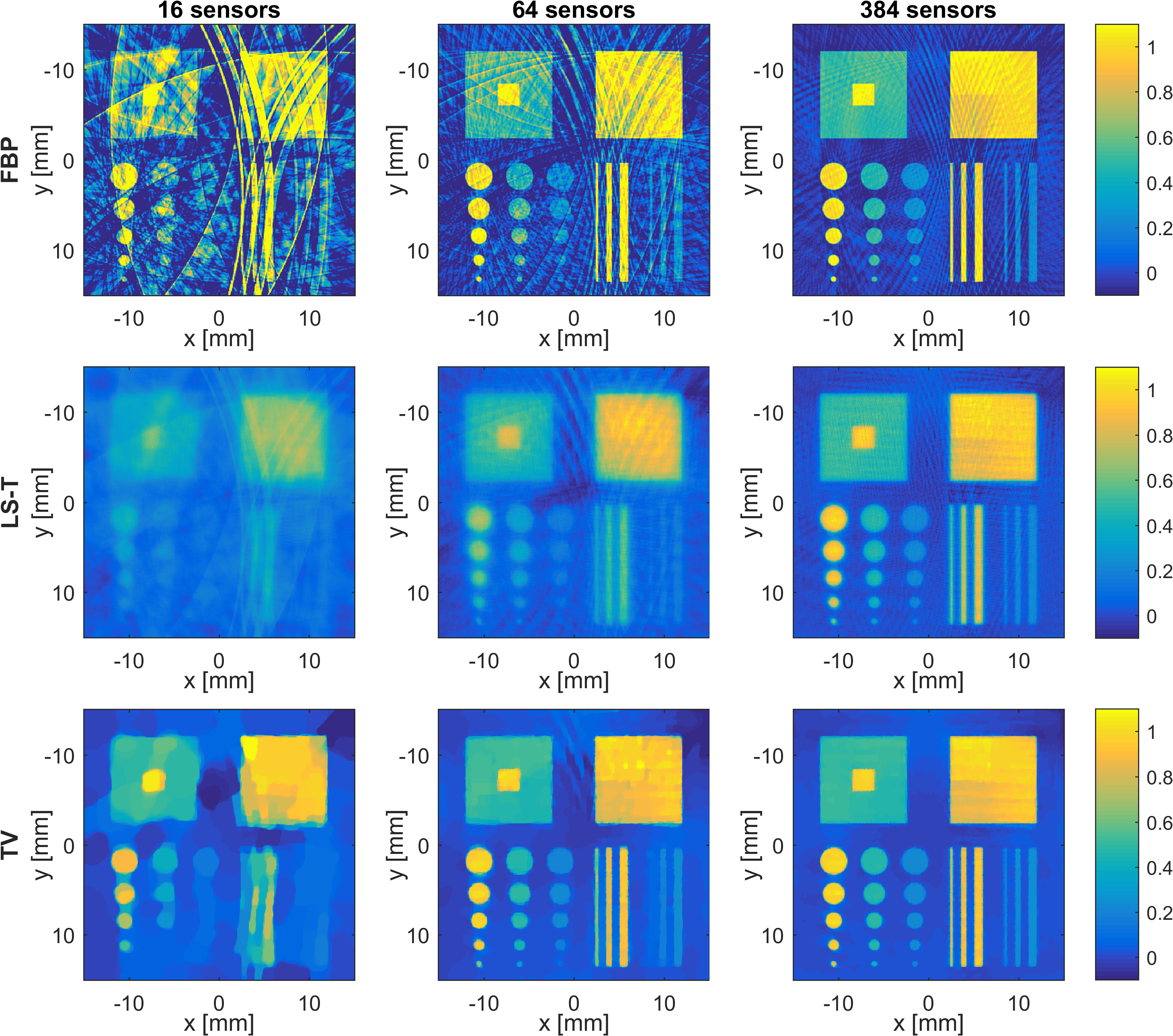}
\caption{Reconstructions of the `piecewise constant' digital phantom in Figure \ref{fig:synth_diverse} with different uniform samplings of the detectors.}
\label{fig:SS_matrix}
\end{figure*}

In Figure \ref{fig:SS_matrix} the results of three reconstruction methods (FBP, LS-T, TV) are shown for coarse sampling (16 detectors), moderate sampling (64 detectors) and fine sampling (384 detectors). All methods give a good reconstruction in case of high sampling, although every method has its own reconstruction bias: TV gives contrast loss in smaller structures, LS-T gives blurred structure edges, whereas FBP gives curved line artefact in the whole image. The TV method is able to keep important features for a longer time as we sample with less detectors: with moderate sampling, it is still possible to detect the minor contrast changes in the upper right square, whereas this is not possible in the other reconstructions. With coarse sampling, the TV reconstruction still gives sharp edges and structures with the right intensity, while this is not the case in the other reconstructions. 

\begin{figure}[!ht]
\centering
\includegraphics[width=\linewidth]{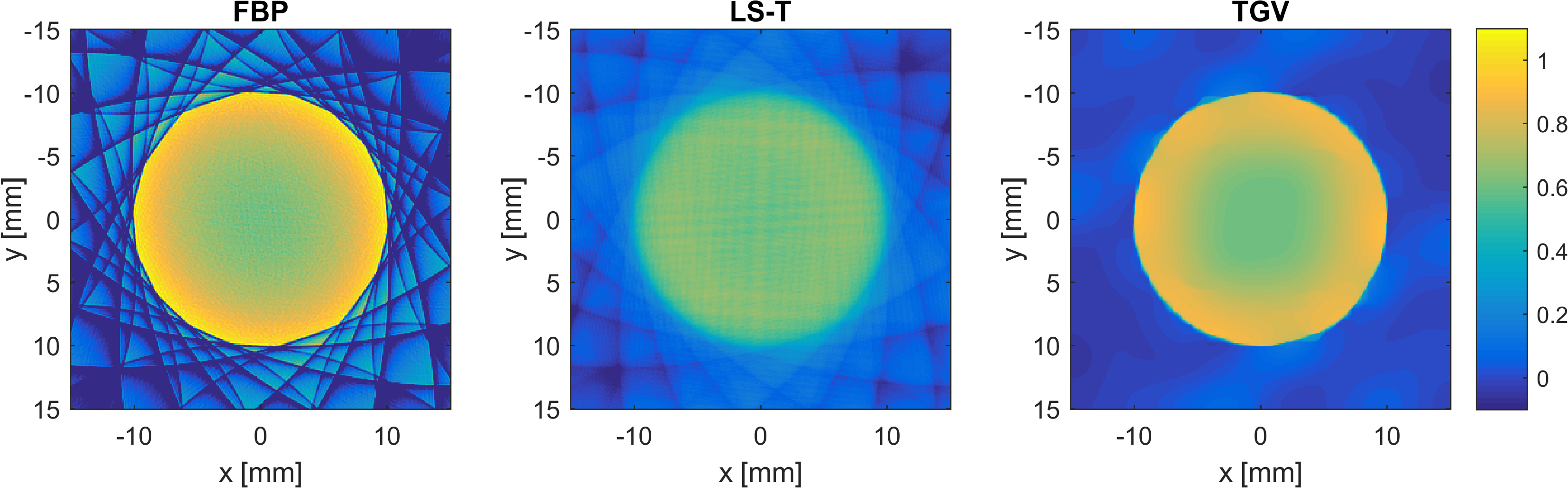}
\caption{Reconstructions of the `smooth disc' digital phantom in Figure \ref{fig:synth_exp} with a uniform sampling of 16 detectors.}
\label{fig:TGV_matrix}
\end{figure}

For the second digital phantom, we again compare with FBP and LS-T. The reconstructions using a uniform sampling of 16 detectors are shown in Figure \ref{fig:TGV_matrix}. The TGV method gives the desired smooth intensity within the disc, while keeping the sharp discontinuities. Moreover, it does not show the curved line artefacts that are visible in the FBP reconstruction.

\begin{figure}[!ht]
\centering
\includegraphics[width=\linewidth]{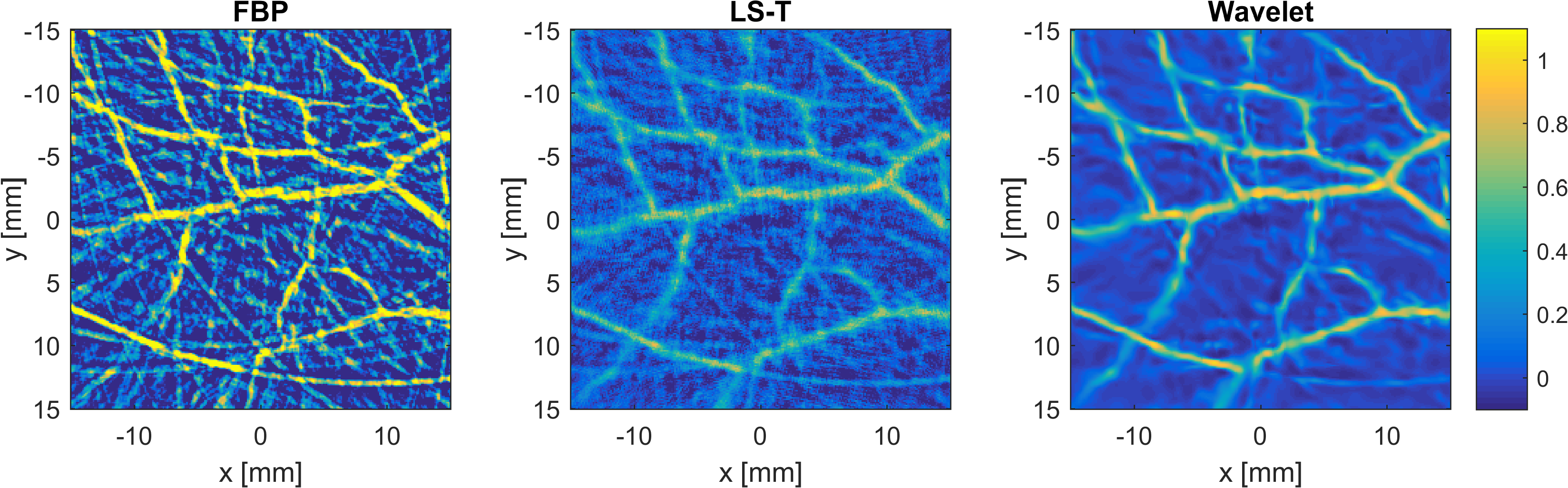}
\caption{Reconstructions of the `vascular structure' digital phantom in Figure \ref{fig:synth_vascular} with a uniform sampling of 32 detectors.}
\label{fig:wavelet_matrix}
\end{figure}

In Figure \ref{fig:wavelet_matrix}, the reconstructions of the third phantom using a uniform sampling of 32 detectors have been shown. FBP performs poorly for this vascular data set, since it is not completely clear which parts of the reconstruction are vascular structures and which are curved line artefacts. The LS-T and wavelet reconstructions show less pronounced curved line artefacts than the FBP reconstruction. The wavelet reconstruction additionally has a higher intensity within the vessels and shows a smoother background. It is striking that regularisation with wavelets is not as effective in removing the curved line artefacts as previous regularisation with TV or TGV, as can be seen around $y=13$ and $x>0$ in Figure \ref{fig:wavelet_matrix}. The probable reason for this is that not only the vascular structure, but also the curved line artefacts can be sparsely represented in the directional wavelet basis. Note that this is only the case in a 2D setting, since in 3D, a backprojection artefact looks like part of the surface of a sphere instead of a curved line. Directional wavelets might be much more effective in removing these 3D artefacts, as long as a suitable transform is used that can sparsely represent anisotropic vascular structures. 

\begin{figure}[!htb]
\centering
\begin{subfigure}{0.47\linewidth}
\includegraphics[width=\linewidth]{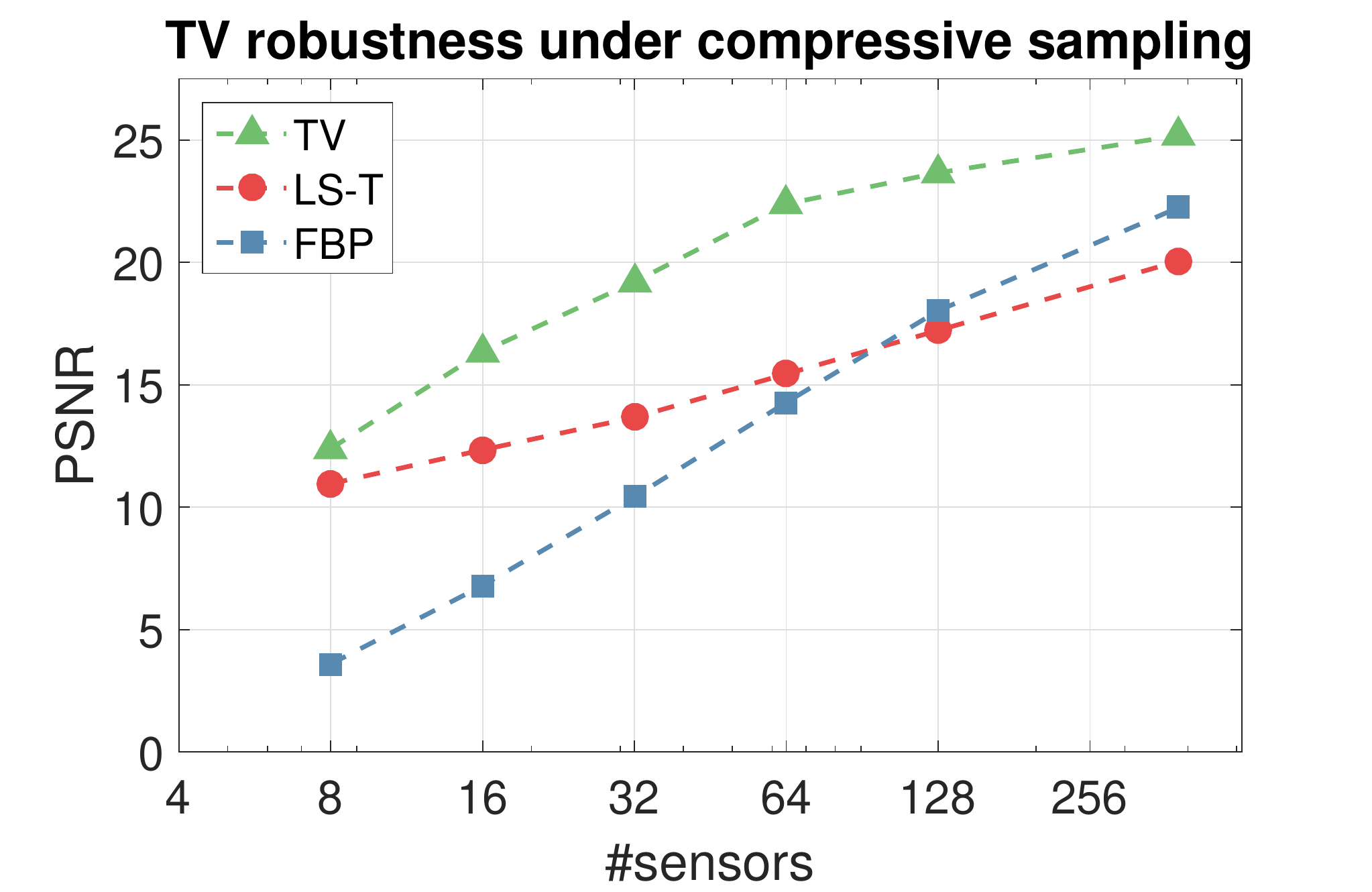}
\caption{} 
\label{fig:SS_diverse}
\end{subfigure}
\begin{subfigure}{0.47\linewidth}
\includegraphics[width=\linewidth]{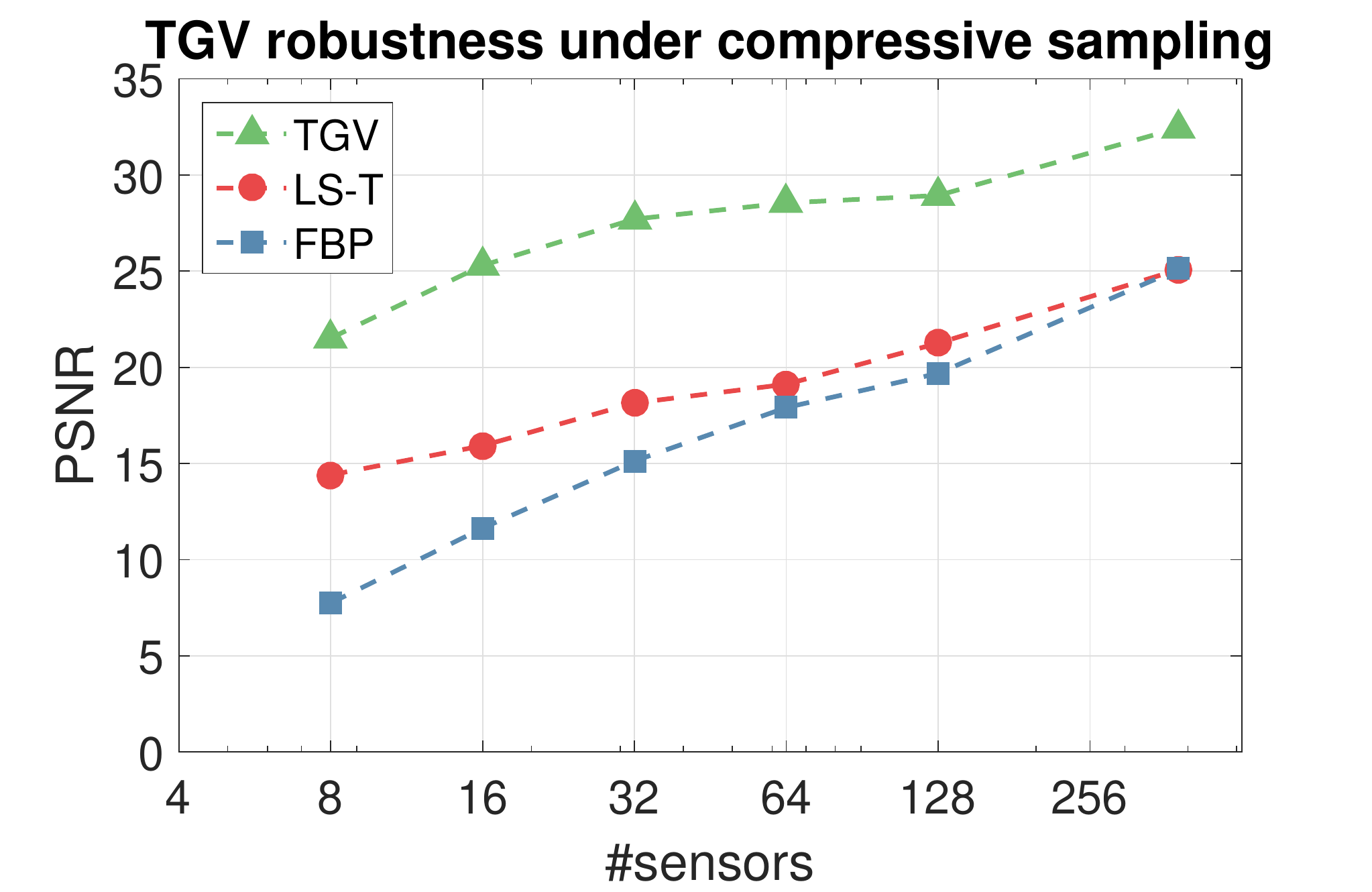}
\caption{~}
\label{fig:SS_exp}
\end{subfigure}
\begin{subfigure}{0.47\linewidth}
\includegraphics[width=\linewidth]{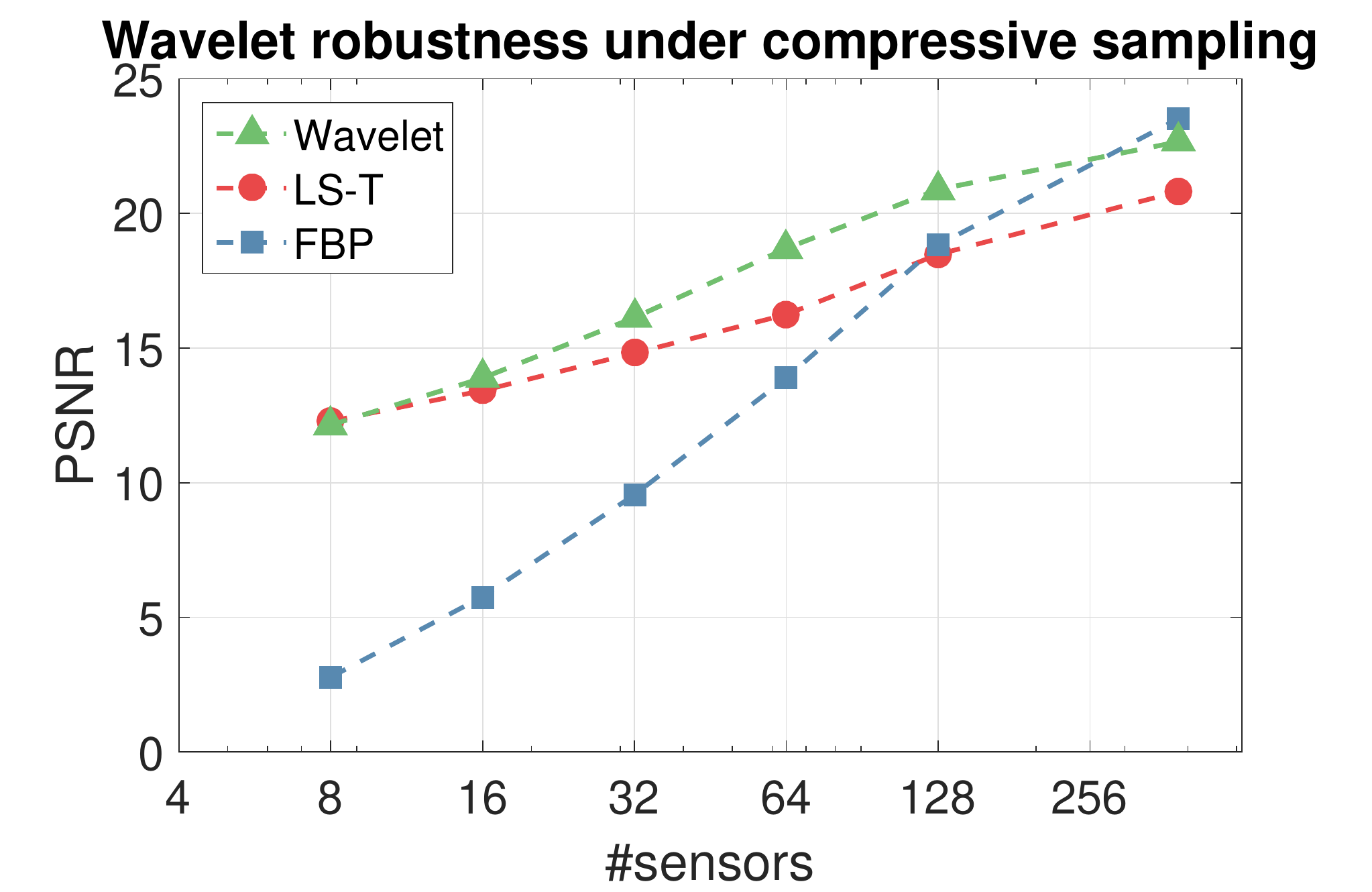}
\caption{~}
\label{fig:SS_vascular}
\end{subfigure}
\caption{PSNR values for different reconstruction methods applied to compressively sampled data. Simulated data from digital phantom in (a) Figure \ref{fig:synth_diverse}; (b) Figure \ref{fig:synth_exp}; (c) Figure \ref{fig:synth_vascular}.}
\label{fig:SS_plot}
\end{figure}

In Figure \ref{fig:SS_plot}, a plot of the PSNR values for different reconstruction methods under compressive sampling is shown. Both the TV and the TGV method perform better than FBP and LS-T. The wavelet method gives a minor improvement under compressive sampling, probably due to the problem that artefacts are too similar to the vascular structure. No noise has been added to the data for this comparison. 

\subsection{Robustness against noise}
\begin{figure}[!htb]
\centering
\begin{subfigure}{0.47\linewidth}
\includegraphics[width=\linewidth]{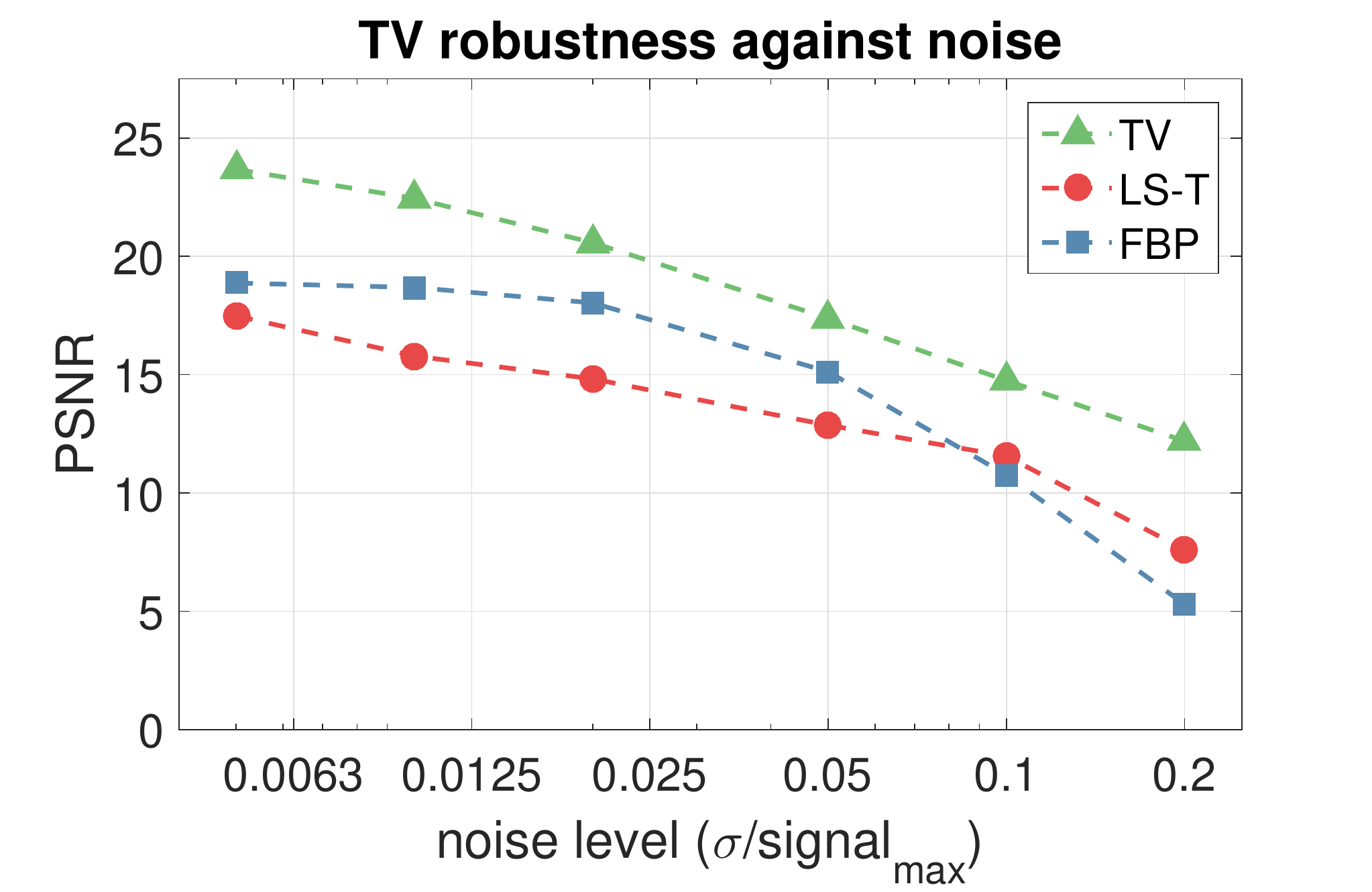}
\caption{} 
\label{fig:NL_diverse}
\end{subfigure}
~~~~
\begin{subfigure}{0.47\linewidth}
\includegraphics[width=\linewidth]{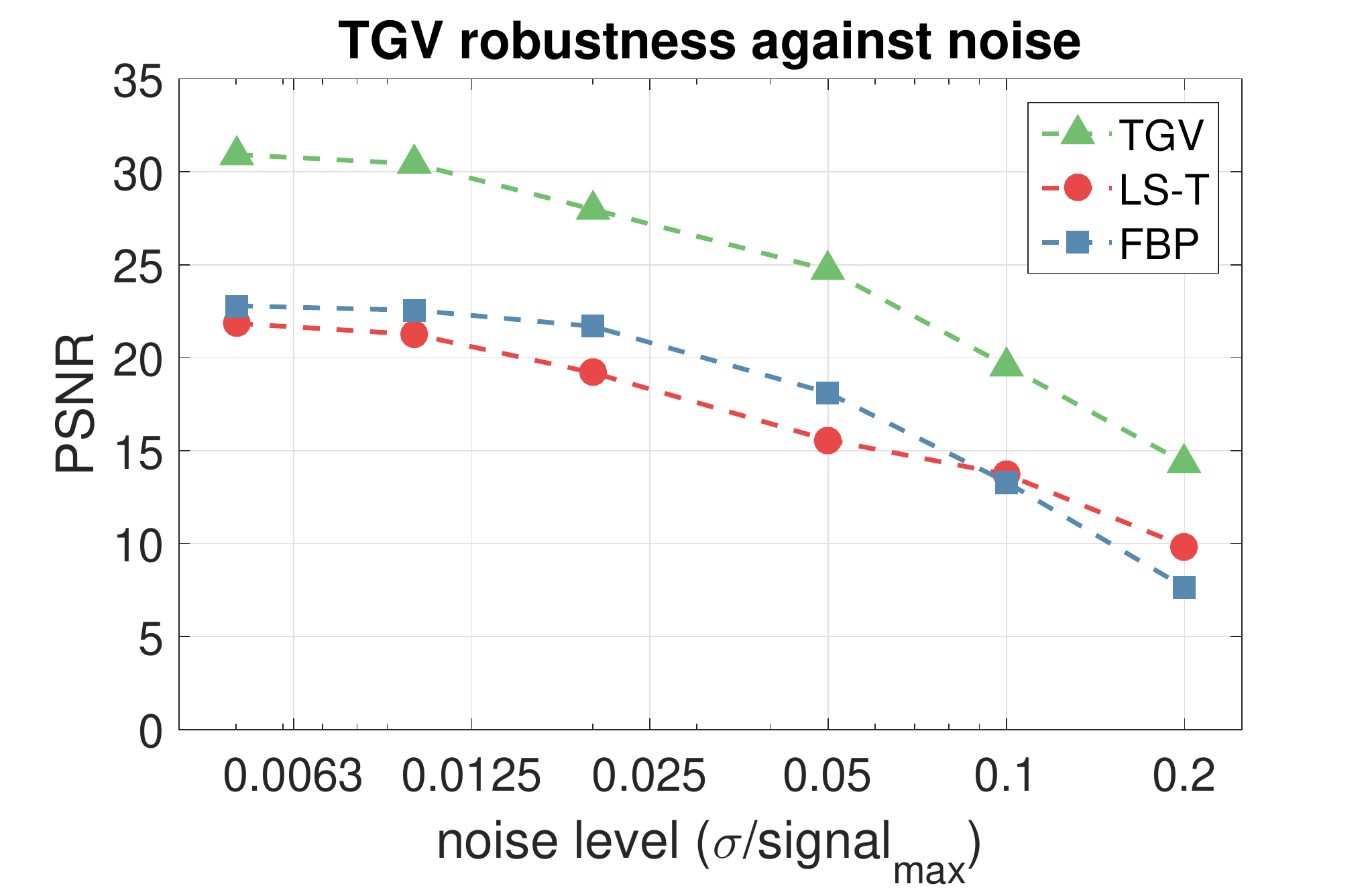}
\caption{~}
\label{fig:NL_exp}
\end{subfigure}
~~~~
\begin{subfigure}{0.47\linewidth}
\includegraphics[width=\linewidth]{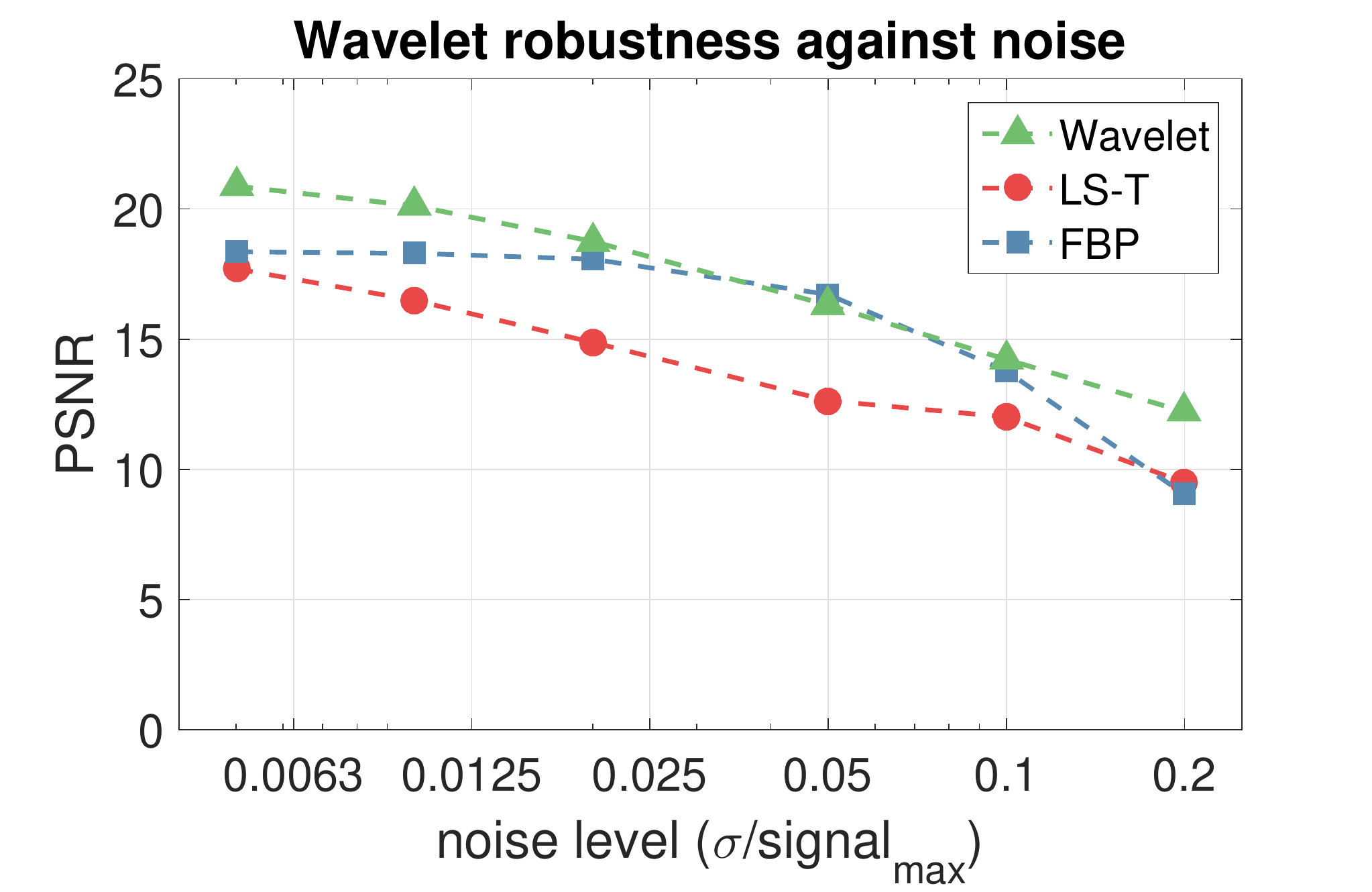}
\caption{~}
\label{fig:NL_vascular}
\end{subfigure}
\caption{PSNR values for different reconstruction methods applied to data with additive Gaussian noise. Simulated data from digital phantom in (a) Figure \ref{fig:synth_diverse}; (b) Figure \ref{fig:synth_exp}; (c) Figure \ref{fig:synth_vascular}.}
\label{fig:NL_all}
\end{figure}

As explained in section \ref{sec:var_meth}, we expect mainly additive Gaussian noise, because of the system electronics and thermal noise from the transducers. After the simulated measured pressure $p(\mathbf{x},t)$ has been obtained, Gaussian noise with zero mean and standard deviation $\sigma$ was added to the measured pressure. In Figure \ref{fig:NL_all}, a comparison between reconstruction methods for the noisy data has been shown. The TV, TGV and wavelet method outperform FBP and LS-T for both low and high noise levels. The results were obtained with a uniform sampling of 192 detectors. 

\section{Experimental results}\label{sec:results_exp}
For the acquisition of experimental data, the cylindrical phantom was scanned in six rotations, giving a total uniform sampling of 384 detectors. We compare the TV, TGV and wavelet reconstruction methods with FBP by using both the full data (384 detectors) and using 16.7\% compressively sampled data (64 detectors) for the reconstructions. The LS-T reconstruction is omitted, since under low noise and middle to high sampling, FBP has shown to outperform LS-T for the synthetic case. The regularisation parameters for the regularised reconstruction have been chosen similar to the parameters found in the synthetic tests, with a small deviation such that the reconstruction is visually best. 

\begin{figure}[!htb]
\centering
\includegraphics[width=0.7\linewidth]{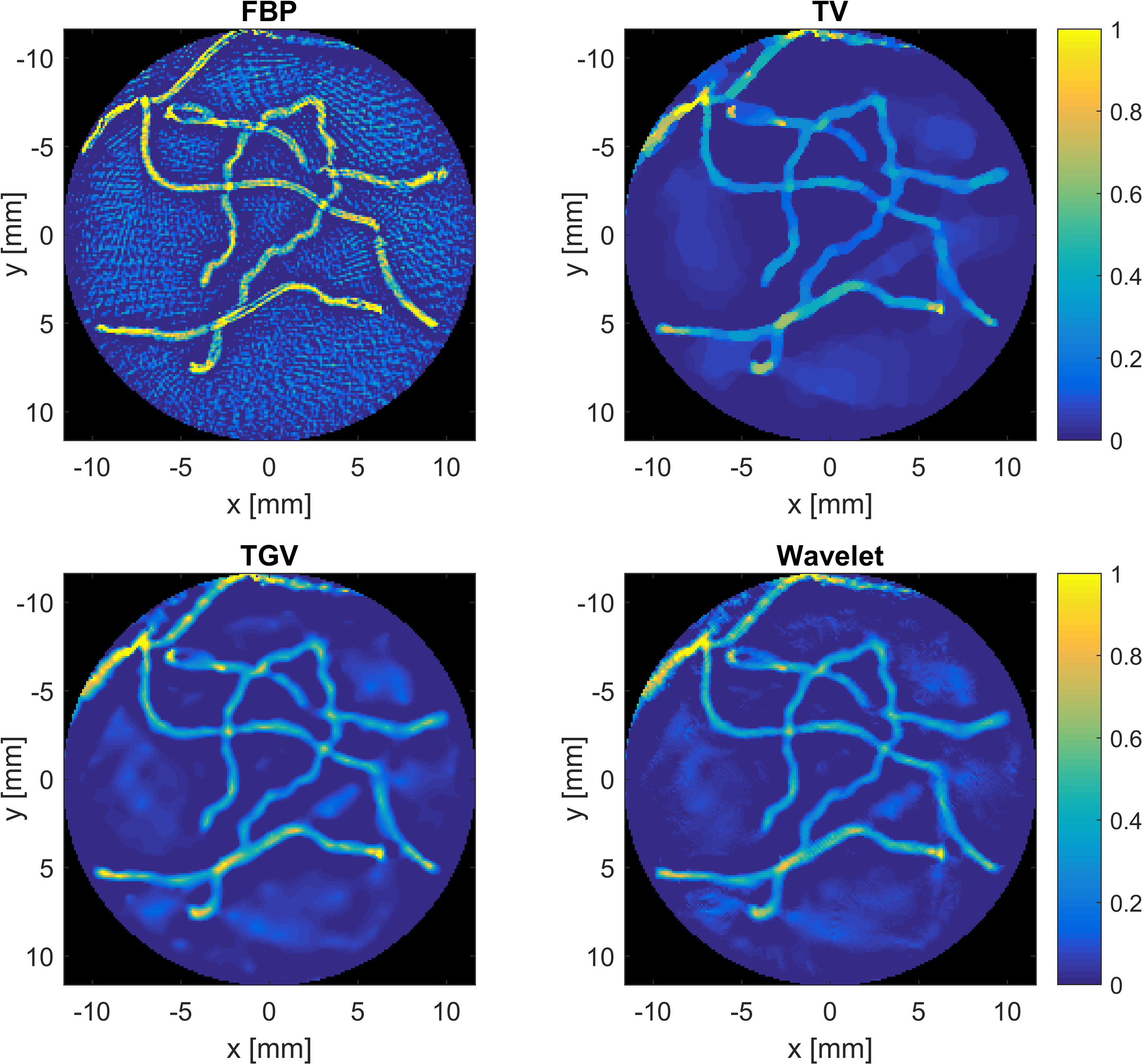}
\caption{Reconstructions of the experimental phantom in Figure \ref{fig:exp} with a uniform sampling of 384 detectors (full reconstruction).}
\label{fig:Exp_vascular_full}
\end{figure}

In Figure \ref{fig:Exp_vascular_full} it can be seen that FBP gives a sharp reconstruction with the expected curved line artefacts in the background. Due to these artefacts, the intensity of the anisotropic structure is not uniform along the structure and undesired low intensity gaps are created. All three regularisation methods give a much smoother reconstruction, both in the background and in the anisotropic structure. Moreover, the small gaps of the FBP reconstruction are filled. The methods give the expected result: a piecewise constant reconstruction for TV, a piecewise linear reconstruction for TGV and an anisotropy enhanced reconstruction for wavelet. Since the dual-tree wavelet approach makes use of smooth elongated structure as it's `building blocks', this reconstruction looks very similar to the TGV reconstruction, although the anisotropic structure shows a slightly sharper boundary in the wavelet reconstruction. 

\begin{figure}[!htb]
\centering
\includegraphics[width=0.7\linewidth]{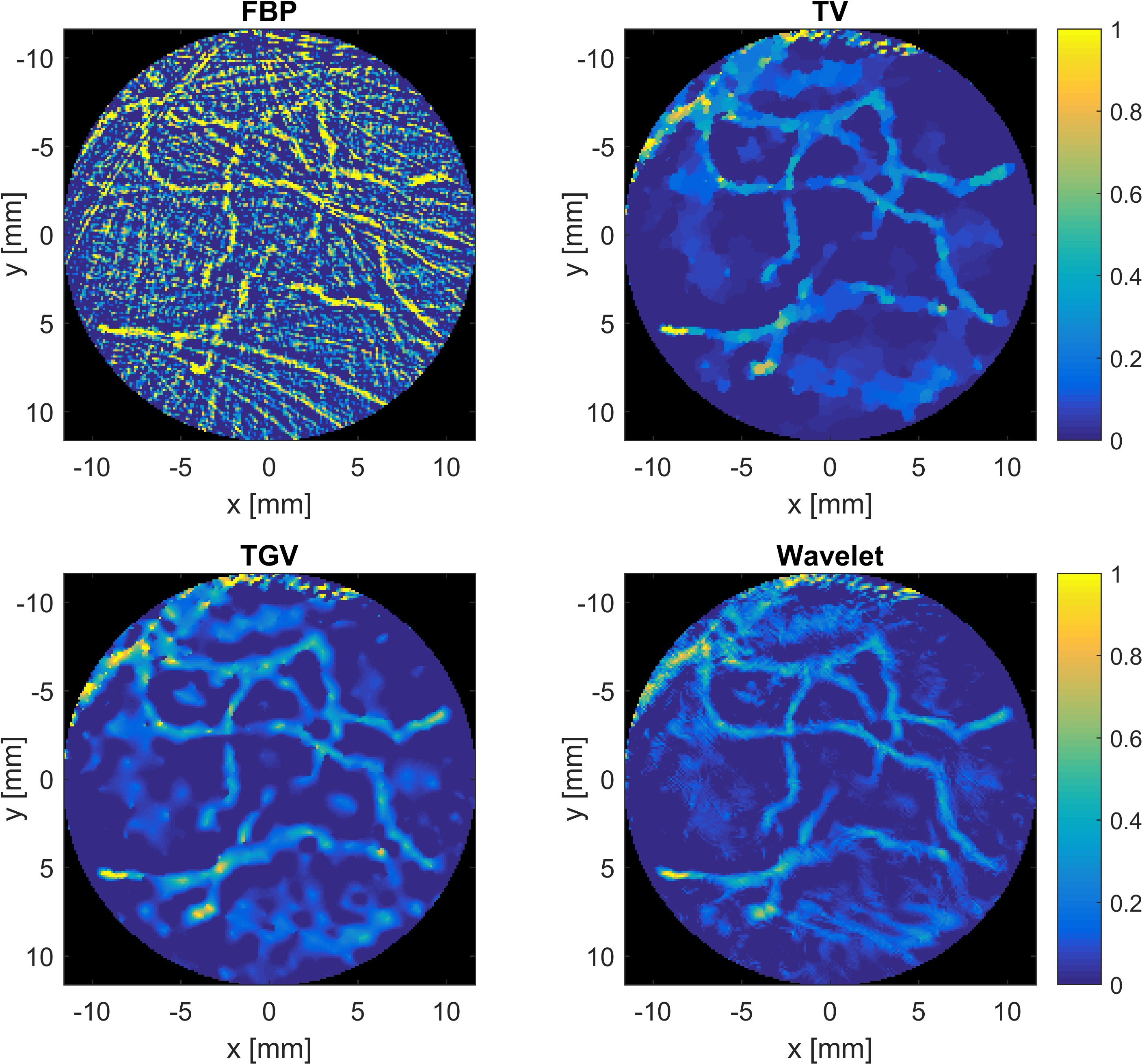}
\caption{Reconstructions of the experimental phantom in Figure \ref{fig:exp} with a uniform sampling of 64 detectors (compressive sampling reconstruction).}
\label{fig:Exp_vascular_CS}
\end{figure}

From the compressive sampling reconstruction (Figure \ref{fig:Exp_vascular_CS}) we can draw similar conclusions: the overall reconstructions are smoother and many gaps that are apparent in the FBP reconstruction are filled in the other reconstructions. If we focus on the region below $y=5$ and between $x=0$ and $x=5$, it is clear that FBP gives strong curved line artefacts, while with the other methods these artefacts are clearly reduced. In these reconstructions, there is much less similarity between the TGV and wavelet results. Due to the use of direction elements, the wavelet reconstruction presents a much better connected structure than the TGV reconstruction. 

\begin{figure}[!htb]
\centering
\begin{subfigure}{0.47\linewidth}
\includegraphics[width=\linewidth]{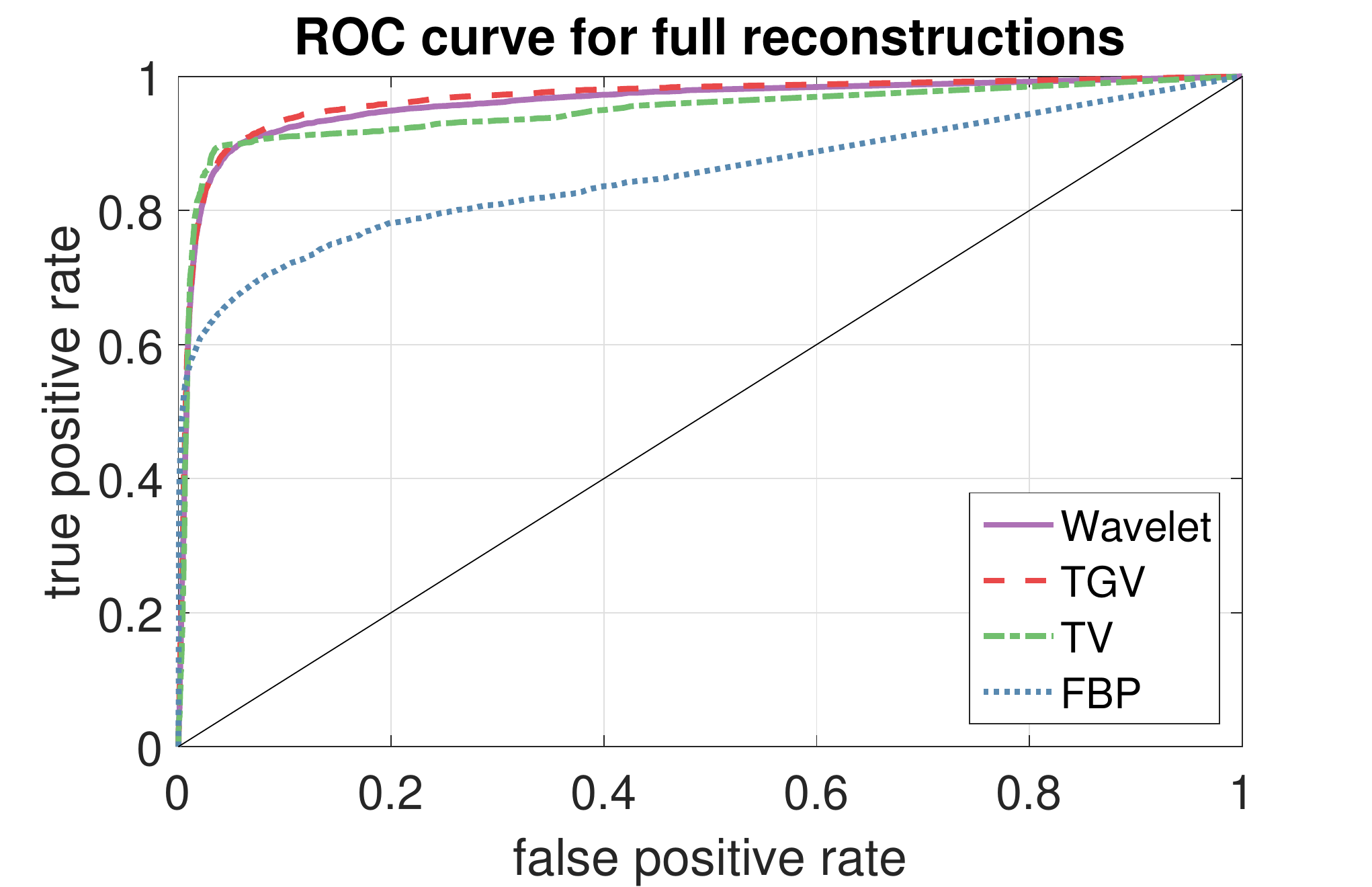}
\caption{} 
\label{fig:ROC_full}
\end{subfigure}
~~~~
\begin{subfigure}{0.47\linewidth}
\includegraphics[width=\linewidth]{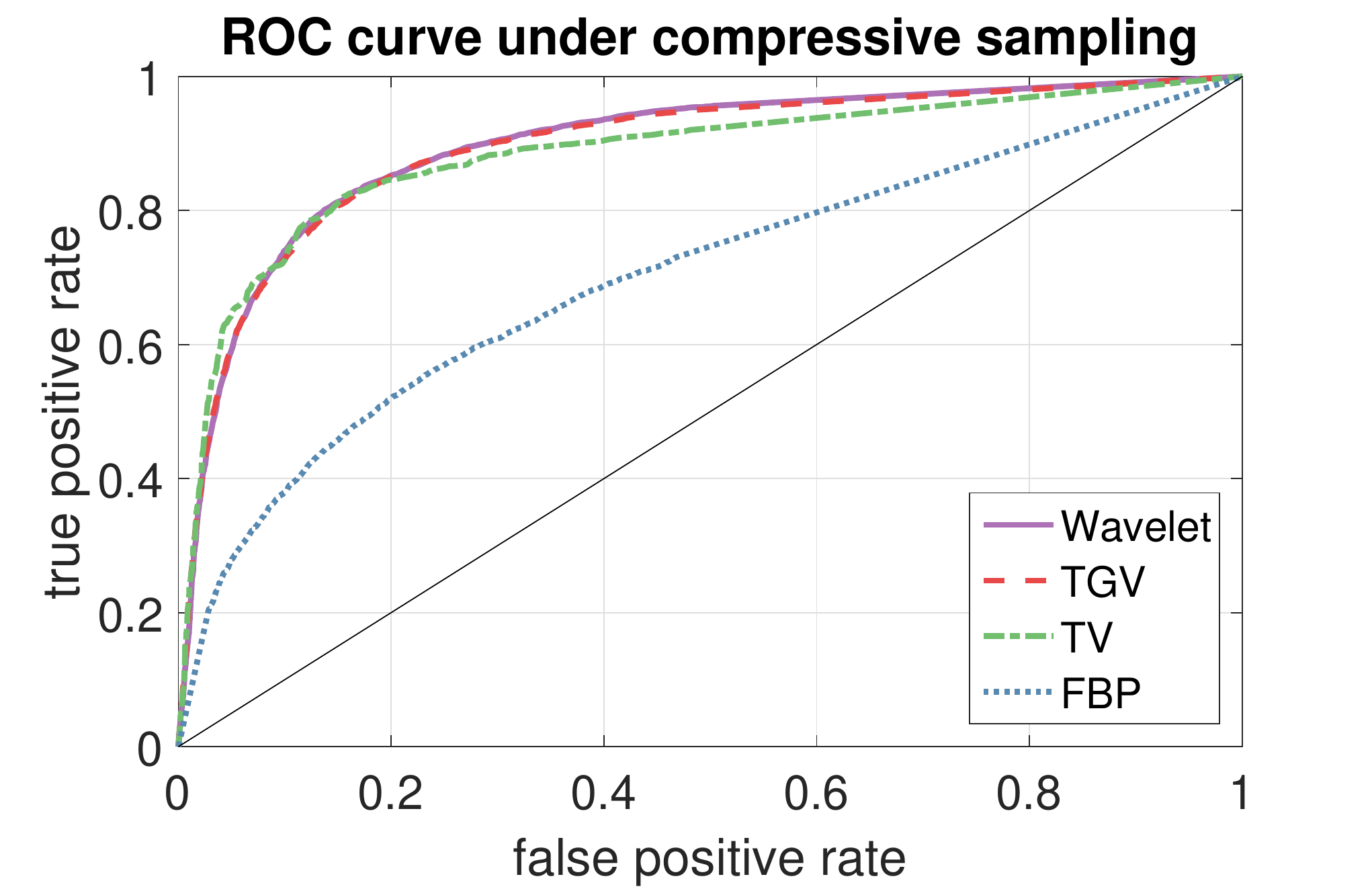}
\caption{~}
\label{fig:ROC_CS}
\end{subfigure}
\caption{ROC curves of different reconstruction methods for full and compressive sampling reconstruction. Curves were obtained by varying the segmentation threshold of the reconstructions and comparing it with the ground truth segmentation of Figure \ref{fig:exp_segmented}.}
\label{fig:ROC_all}
\end{figure}

As can be seen in Figure \ref{fig:ROC_full}, the ROC curves of the regularised reconstructions show an almost ideal behaviour. The ROC curve of the FBP reconstruction lies much lower, which shows that this method is a poor choice when it comes to thresholding capability. In Figure \ref{fig:ROC_CS}, we see a similar relation between the reconstruction methods, although all curves lie lower than in Figure \ref{fig:ROC_full}. It is interesting to see that all regularised reconstruction methods give a similar ROC curve. TGV and wavelet almost completely overlap, while TV is a bit higher in the steep part, but continues slightly lower in the second part of the graph. This tells us that either method will help us to correctly segment the photoacoustic reconstruction for anisotropic vascular structures.

\section{Summary and outlook}\label{sec:discussion}
In this work we have derived a general variational framework for regularised PAT reconstruction where no specific measurement geometry or physics modelling is assumed. The framework can be applied to a PAT reconstruction problem in a 2- or 3-dimensional setting. The primal-dual implementation is a modular one and the optimisation parameters are chosen such that the algorithm is efficient. Specific choices for the regulariser made in this paper are TV, second order TGV and the complex dual-tree wavelet transform. These regularisers can respectively deal with the prior assumptions of either sharp discontinuities in combination with a homogeneous and heterogeneous fluence rate, or anisotropic structures. The methods have been compared on both synthetic and experimental data with two widely used direct reconstruction methods.

We have shown that the modularity of the variational framework enables easy exchange of regularisers, which arise from different prior assumptions on the reconstruction. Using this framework, one can choose the terms within as desired, without having to change the structure of the algorithm.

Furthermore it was shown that our TV and TGV methods perform better than direct reconstruction methods: they are better able to handle both low and high noise levels and give better reconstructions under uniform compressive sampling. In a clinical setup, it might not be preferred to change the detector locations during a full scan, because of limited view or movement of the tissue under consideration. For this reason, the use of TV and TGV methods are promising. In this work, reconstructions were made using uniform compressive sampling. Future efforts could lie in finding the best sampling strategy for various data sets and regularisers.

Finally, research could be focussed on the reconstruction of anisotropic structures such as blood vessels and the difference between this reconstruction in 2D and 3D. It is unclear if the similarity between backprojection artefacts and vascular structures is preventing better reconstructions in the 2D case and if this problem is absent in the 3D case. For this reason, it is useful to develop or learn a dictionary that allows sparse reconstruction of vascular structures in 3D.

\section*{Acknowledgements}
The authors thank Maura Dantuma for the acquisition of the experimental data. Marinus J. Lagerwerf acknowledges financial support from the Netherlands Organization for Scientific Research (NWO), project  639.073.506. The authors acknowledge Stichting Achmea Gezondheidszorg for funding in project Z620, and the Pioneers in Healthcare Innovation (PIHC) fund 2014 for funding in project RAPACT.

\small{
\bibliographystyle{abbrv} % alternative abbrv, alpha, amsalpha, amsplain 
\bibliography{refs}
}
\end{document}

%% file: packages.tex
\usepackage{subfiles}
\usepackage{algpseudocode}
\usepackage{arydshln}

\usepackage{fullpage} %http://en.wikibooks.org/wiki/LaTeX/Page_Layout
\usepackage[affil-it]{authblk}%

\usepackage[pagebackref,pdftex,bookmarks,colorlinks,breaklinks]{hyperref}
\hypersetup{%
		linkcolor=magenta,      %has effects if colorlinks=true
		citecolor=blue,				
		filecolor=dullmagenta,
		urlcolor=magenta,
}

\usepackage[T1]{fontenc}    % Zeichensatzkodierung T1
\usepackage[utf8]{inputenc} % Definition des Eingabeencodings, also wie Umlaute umgesetzt werden
\usepackage{lmodern}
\usepackage{ae}             % Verwendung von AE-Fonts
\usepackage{amsmath,amstext,amsfonts,amssymb,amsthm,latexsym}
\usepackage{todonotes}

\usepackage{graphicx}         % for using includegraphics
\usepackage{epsfig}           % for using epsfig <-- uses includegraphics
\usepackage{subcaption}
\usepackage{algorithm}
\usepackage{algpseudocode}
\usepackage{bbm}

%% file: commands.tex
\newcommand{\R}{\mathbb{R}}

\newcommand{\F}{\mathcal{F}}
\newcommand{\C}{\mathcal{C}}

\newcommand{\BV}{\text{BV}(\Omega)}

\newcommand{\diff}{\textnormal{\,d}}

\renewcommand{\div}{\text{div}}

\newcommand{\argmin}{\operatornamewithlimits{argmin}}

\newcommand{\prox}{\text{prox}}

\newcommand{\norm}[1]{\| #1 \|}

\newcommand{\bigo}{\mathcal{O}}

\newcommand{\bbm}{\begin{pmatrix}[l]}
\newcommand{\ebm}{\end{pmatrix}}

\newcommand{\grad}{\nabla}%
\newcommand{\TV}{\text{TV}}
\newcommand{\TGV}{\text{TGV}}